\documentclass[final]{siamltex}
\usepackage{amsfonts,amsmath,mathrsfs,bm}
\usepackage{amssymb}
\usepackage[notref,notcite]{showkeys}
\usepackage{hyperref}
\usepackage{graphicx}
\usepackage{wrapfig}
\usepackage{subfig}
\usepackage{float}
\usepackage{bbm}
\usepackage{algorithmic, algorithm}
\usepackage{gensymb,color,soul}
\usepackage{tikz}
\usepackage{tikz-3dplot}
\usepackage{mathtools}
\mathtoolsset{showonlyrefs}

\newcommand{\N}{\mathbb{N}}
\newcommand{\R}{\mathbb{R}}
\newcommand{\C}{\mathbb{C}}

\newtheorem{remark}[theorem]{Remark}

\title{Local electrical impedance tomography via projections}

\author{
A.~J\"a\"askel\"ainen\footnotemark[2]
\and A.~Vavilov\footnotemark[2]
\and J.~Toivanen\footnotemark[3]
\and A.~Hänninen\footnotemark[3]
\and V.~Kolehmainen\footnotemark[3]
\and N.~Hyv\"onen\footnotemark[2]
}

\begin{document}
\maketitle

\renewcommand{\thefootnote}{\fnsymbol{footnote}}
\footnotetext[2]{Aalto University, Department of Mathematics and Systems Analysis, P.O.~Box 11100, FI-00076 Aalto, Finland (altti.jaaskelainen@aalto.fi, anton.vavilov@aalto.fi, nuutti.hyvonen@aalto.fi). This work was supported by the Jane and Aatos Erkko Foundation and the Academy of Finland (decision 348503, 353081, 359434, 359181).}
\footnotetext[3]{
University of Eastern Finland, Department of Technical Physics, Kuopio Campus, P.O. Box 1627, FI-70211 Kuopio, Finland (jussi.toivanen@uef.fi, asko.hanninen@uef.fi, ville.kolehmainen@uef.fi). This work was supported by the Jane and Aatos Erkko Foundation and the Academy of Finland (decision 353084, 359433, 358944).}

\begin{abstract}
This paper introduces a method for approximately eliminating the effect that conductivity changes outside the region of interest have in electrical impedance tomography, allowing to form a local reconstruction in the region of interest only. The method considers the Jacobian matrix of the forward map,~i.e.,~of the map that sends the discretized conductivity to the electrode measurements, at an initial guess for the conductivity. The Jacobian matrix is divided columnwise into two parts: one corresponding to the region of interest and a nuisance Jacobian corresponding to the rest of the domain. The leading idea is to project both the electrode measurements and the forward map onto the orthogonal complement of the span of a number of left-hand singular vectors for a suitably weighted nuisance Jacobian. The weighting can,~e.g.,~account for the element sizes in a finite element discretization or to prior information on the conductivity outside the region of interest. The inverse problem is then solved by considering the projected relation between the measurements and the forward map, only reconstructing the conductivity in the region of interest. The functionality of the method is demonstrated by applying a reconstruction algorithm that combines lagged diffusivity iteration and total variation regularization to experimental data. In particular, data from a head-shaped water tank is considered, with the conductivity change in the region of interest mimicking growth of a hemorrhagic stroke and the changes outside the region of interest imitating physiological variations in the conductivity of the scalp.  
\end{abstract}

\renewcommand{\thefootnote}{\arabic{footnote}}



\pagestyle{myheadings}
\thispagestyle{plain}
\markboth{A.~J\"A\"ASKEL\"AINEN ET AL.}{LOCAL ELECTRICAL IMPEDANCE TOMOGRAPHY VIA PROJECTIONS}

\section{Introduction}
\label{sec:introduction}

Electrical impedance tomography (EIT) is an imaging method for studying the internal conductivity of an object. The method is based on driving electrical current through the object via electrodes placed on its surface and measuring the resulting potential differences between the electrodes. Given the driven currents as well as the measured potentials, an estimate for the internal conductivity of the object can be formed by solving a nonlinear inverse elliptic boundary value problem. For more information on EIT, see~\cite{Borcea02,Cheney99,Uhlmann09}.

In this work, the EIT measurement process is modeled by a smoothened version of the complete electrode model (CEM), which has been found to correspond well with real-world measurements \cite{Cheng89}. The CEM models electrodes as subsets of the boundary of the examined object, and its basic version assumes a constant contact conductivity on each electrode surface~\cite{Somersalo92}. The smoothened CEM (SCEM) differs from this by assuming that the contact conductivity is a smooth function that vanish at the electrode boundaries. The advantage of SCEM is that it leads to higher regularity for the forward solution, which has potential to induce better convergence for numerical methods for solving the forward problem of EIT~\cite{Hyvonen17b}.

A potential application of EIT is head imaging, particularly for the purpose of on-line monitoring of a stroke patient in intensive care. While a computerized tomography (CT) measurement is a standard tool for stroke imaging and diagnosis, it can not be used
for on-line monitoring. EIT is potentially a well suited solution for on-line detection of temporal changes in the brain, such as growth of intracerebral hemorrhagic stroke or secondary bleeding in ischemic stroke patients, in a bedside monitoring setting~\cite{Toivanen21,Toivanen24}. This is achieved by utilizing (linear) difference imaging to reconstruct the difference in conductivity between two measurements taken at different times; see, e.g.,~\cite{Adler09,Barber84}. The conductivity inside a head can, however, also change due to other reasons that are not indicative of abnormalities. These changes can have a significant effect on the measurements and thus hinder one's ability to accurately reconstruct developments of the stroke. An example of this is blood circulation in the skin, which can have a major impact on the measurements due to its close proximity to the electrodes.

In this work, we utilize projections to partially eliminate effects of conductivity changes in regions we are not interested in, allowing to (more accurately) reconstruct the conductivity in the {\em region of interest} (ROI) only. This method is an extension of our previous approach  to reduce the effect of unknown electrode contacts via projections in \cite{Jaaskelainen25}, where the forward problem was projected onto the orthogonal complement of the range of the Jacobian matrix of the potential measurements with respect to the electrode contacts. When computed with respect to a discretized conductivity over a {\em region of non-interest} (RONI) rather than electrode contacts, the range of such a {\em nuisance Jacobian} has a much higher dimension, often coinciding with the whole image space. As a consequence, a projection onto the orthogonal complement of the range of the Jacobian removes (almost) all information in the measurements, which renders our previous approach inapplicable as such. To deal with this, we use only a partial projection constructed by choosing a subspace of the range of the nuisance Jacobian with the highest expected influence on the measurements, given some prior information,~cf.~\cite{Calvetti25}. Such a subspace can be constructed by computing a number of left singular vectors for a suitably weighted nuisance Jacobian, with the weight matrix accounting for,~e.g.,~the sizes of elements in a finite element (FE) discretization or prior information on the conductivity inside the RONI.

Approximate marginalization of conductivity changes in the RONI and closely related problems have previously been tackled by employing the approximation error method~\cite{Kaipio05}; see,~e.g.,~\cite{Hadwin14,kaipio2013,Nissinen09}. The method accounts for the effect of conductivity changes in the RONI by incorporating additive Gaussian approximation error noise into the measurement model. Before the measurements are available, the second order statistics for the approximation error noise can be estimated via simulations based on random draws from the priors for the conductivity and other free parameters. Subsequently, the approximation error noise can be coupled with standard Gaussian measurement noise in Bayesian inversion. Compared to the approximation error approach, a strength of our projection method is that it does not require extra random draws from the priors or the computation of related forward solutions.

We study the effectiveness of our projection method using real-world measurements from water tanks. Measurements are taken from two cylindrical tanks, where our goal is to test the method in simple geometries, and a head-shaped tank. The head tank contains a resistive 3D-printed perforated skull, and we use projections to reduce the effect of conductivity perturbations in the layer exterior to the skull, representing the scalp.

The measurements are analyzed by computing reconstructions of the internal conductivity distributions both with and without projections. Forming a reconstruction constitutes an inverse problem that we solve by resorting to a Bayesian approach after linearizing the forward problem defined by the SCEM at an initial guess for the internal and contact conductivities. We assume a smoothened total variation (TV) prior distribution for the perturbation of the internal conductivity in the ROI (cf.~\cite{Rudin92}) and compute reconstructions using an iterative lagged diffusivity algorithm~\cite{Vogel96,Harhanen15}.  For the simplest example in a disk-shaped water tank, the reference measurements required by the difference imaging approach are computed based on prior information on the internal and contact conductivities. For the more involved experiments, the reference data are obtained via actual measurements with no embedded inclusions present. In some reconstructions, we also utilize projections onto the orthogonal complements of the ranges of Jacobian matrices with respect to the electrode contacts or positions in order to alleviate geometric mismodeling or reduce the effect of inaccurate initial guesses for the contacts, cf.~\cite{Jaaskelainen25}.

This article is organized as follows. Section~\ref{sec:CEM} introduces the SCEM and explains how derivatives of the measured potentials with respect to conductivity perturbations can be computed. In Section~\ref{sec:comp}, we describe the computational models as well as the experimental setup used for performing the measurements. Section~\ref{sec:data} considers our projection method and how it can be incorporated into reconstructions. The algorithms for computing reconstructions are detailed in Section~\ref{sec:algorithms}, while Section~\ref{sec:experiments} contains descriptions and results of our experiments. Finally, our concluding remarks can be found in Section~\ref{sec:conclusion}.

\section{Complete electrode model}
\label{sec:CEM}

The forward problem of finding the measured voltages for given input currents is modeled by a version of the CEM~\cite{Cheng89}. Specifically, the model takes a current pattern $I \in \C^M_\diamond$, which consists of the net currents that are driven through each of $M \in \N \setminus \{ 1 \}$ electrodes, and determines the potentials on each electrode $U \in \C^M_\diamond$. Here $\C^M_\diamond$ is a zero-mean vector space defined as
\begin{equation}
\label{eq:zero-mean}
\C^M_\diamond = \Big\{V\in\C^M\,\Big|\, \sum_{m=1}^M V_m = 0\Big\}.
\end{equation}
While the zero-mean property of the current pattern $I$ is a result of the conservation of charge, the potentials $U$ obtain this property by our choice of the ground potential level. 

The CEM models the imaged object as a bounded Lipschitz domain $\Omega \subset \R^3$, with the electrodes represented as open connected subsets $E_1, \dots ,E_M$ of the boundary $\partial \Omega$ and their union denoted by $E$. It is assumed that the closures of the electrodes are mutually disjoint. With these definitions, the CEM defines the relationship between the currents $I$ and potentials $U$ via an elliptic boundary value problem
\begin{equation}
\label{eq:cemeqs}
\begin{array}{ll}
\displaystyle{\nabla \cdot(\sigma\nabla u) = 0 \qquad}  &{\rm in}\;\; \Omega, \\[6pt]
{\displaystyle {\nu\cdot\sigma\nabla u} = \zeta (U - u) } \qquad &{\rm on}\;\; \partial \Omega, \\[2pt]
{\displaystyle \int_{E_m}\nu\cdot\sigma\nabla u\,{\rm d}S} = I_m, \qquad & m=1,\ldots,M. \\[4pt]
\end{array}
\end{equation}
Here, $U$ is identified with a piecewise constant function on the electrodes, $\nu$ is the exterior unit normal vector of the boundary $\partial \Omega$, and $\sigma \in L^\infty_+(\Omega)$ denotes the internal conductivity of the object, which is assumed to be isotropic, with
\begin{equation}
\label{eq:sigma}
L^\infty_+(\Omega) := \{ \varsigma \in L^\infty(\Omega) \ | \ {\rm ess} \inf {\rm Re} (\varsigma) > 0 \}.
\end{equation}

In \eqref{eq:cemeqs}, $\zeta$ describes the contact conductivity between the surfaces of the electrodes and the object. In the traditional CEM, $\zeta$ is assumed to take a constant value on each electrode $E_m$, and its contribution is typically modeled via contact impedance instead of contact conductance. However, for our experiments we instead use the SCEM~\cite{Hyvonen17b}, in which the contact conductivity $\zeta$ is a smooth function in
\begin{equation}
  \label{eq:zeta}
\mathcal{Z} := \big\{ \xi \in L^\infty(E) \  \big| \   {\rm Re}\, \xi \geq 0 \   {\rm and} \  {\rm ess} \sup \! \big( {\rm Re} ( \xi |_{E_m} ) \big) > 0 \  {\rm for} \ {\rm all} \ m= 1, \dots, M \big\}.
\end{equation}
The definition of $\zeta$ is extended onto the entire boundary $\partial \Omega$ via zero-continuation.

For given $\sigma \in L^\infty_+(\Omega)$, $\zeta \in \mathcal{Z}$ and $I \in \C^M_\diamond$, the forward problem \eqref{eq:cemeqs} defines the interior-electrode potential pair uniquely up to a common ground potential level. In more precise mathematical terms, \eqref{eq:cemeqs} has a unique weak solution $(u,U) \in H^1(\Omega) \oplus \C^M_\diamond$, where we have respected the aforementioned choice for the ground level that forces the electrode potentials to live in $\C^M_\diamond$. For more details on the SCEM, including analysis on the higher regularity of $u$, consult \cite{Hyvonen17b}. 

In order to obtain more information about the internal conductivity of the imaged object, EIT measurements are typically taken using several different current patterns. As the function $U(\sigma,\zeta; I) = [ U_1(\sigma,\zeta; I^{(j)}), \dots, U_M(\sigma,\zeta; I^{(j)})]^\top$, defined by \eqref{eq:cemeqs}, is linear with respect to $I$, $M-1$ linearly independent current patterns are enough to gather all available information on a specific measurement setup (up to the measurement noise). The function representing a full set of measurements for $L \in \N$ input current patterns is denoted as
\begin{equation}
\label{eq:forward_map}
\mathcal{U}\big(\sigma,\zeta; I^{(1)}, \dots, I^{(L)}\big) 
= \big[U(\sigma,\zeta; I^{(1)})^{\top}, \dots, U(\sigma,\zeta; I^{(L)})^{\top}\big]^{\top}
\in \C^{ML}.
\end{equation}
In our experiments, $L$ does not equal $M-1$ due to all electrodes not being capable of feeding currents or as a means to improve the signal-to-noise ratio.

In addition to solving the forward problem of evaluating $U(\sigma,\zeta; I)$, its derivatives with respect to the degrees of freedom in the chosen parametrization for the conductivity $\sigma$ also need to be computed in order to form reconstructions and introduce the employed projections. The forward solution of the SCEM is known to be Fr\'echet differentiable with respect to $\sigma$ \cite{Darde21,Hyvonen17b}, and the derivative in the direction $\eta \in L^\infty(\Omega)$ can be efficiently computed with the sampling formula
\begin{equation}
\label{eq:sderiv}
D_\sigma U(\sigma; \eta) \cdot \tilde{I} = - \int_{\Omega} \eta \nabla u \cdot \nabla \tilde{u} \, {\rm d} x,
\end{equation}
where $u$ and $\tilde{u}$ are, respectively, the first parts of the solutions to \eqref{eq:cemeqs} for the current patterns $I$ and $\tilde{I}$. 
Similar sampling formulas also exist for the derivatives with respect to the contact conductivity $\zeta$ and electrode locations~\cite{Darde12,Darde21,Hyvonen14}, and they are employed when computing Jacobian matrices and associated projections with respect to electrode contacts and positions in our numerical experiments.

\section{Computational model and experimental setup}
\label{sec:comp}
Our EIT measurements are performed with low-frequency alternating currents, which naturally leads to a model with complex-valued conductivities, current patterns and electric potentials~\cite{Borcea02} as described in the preceding section. However, we only measure the amplitudes of electrode currents and voltages and utilize them as inputs and outputs for a forward model characterized by real-valued internal and contact conductivities. This approximation has been successfully employed in many experimental EIT studies. In consequence, the multiplier field $\C$ is replaced by the real numbers $\R$ in the rest of this paper.

\subsection{Finite element discretization}
We use a FE discretization of the domain~$\Omega$ for the numerical computations performed in our experiments.  The discretized domain is characterized by a tetrahedral mesh, on which we approximate conductivity and potential distributions using FE basis functions $\varphi_j$ that are associated with the vertices of the mesh numbered as $j = 1,\dots, n$. The basis functions are piecewise linear ``hat'' functions each of which takes the value 1 at a single vertex and the value 0 at all other vertices. An approximation for the internal conductivity of $\Omega$ can then be formed as
\begin{equation}
\label{eq:discr_sigma}
\sigma = \sum_{j=1}^n \sigma_j \varphi_j \, ,
\end{equation}
which allows to abuse the notation by identifying the conductivity with a vector $\sigma = [\sigma_1,\dots,\sigma_n]^\top \in \R^n$. In a similar fashion, we use the FE basis to discretize the internal potential $u$. In the following sections discussing the computational implementation or numerical results, we use $\sigma$ and $u$ to refer to these piecewise linear FE approximations. 

It is assumed that the ROI and RONI cover the whole discretized domain, and the respective nodes of the FE mesh can formally be numbered so that the indices $1, \dots, n_{\rm b}$ correspond to the ROI and $n_{\rm b} + 1, \dots, n =n_{\rm b} + n_{\rm s}$ correspond to the RONI. Here, the subindices b and s refer to ``brain'' and ``skin'' in anticipation of the application to the head-shaped water tank. In particular, we decompose the conductivity vector into two parts as $\sigma^\top = [ \sigma_{\rm b}^\top, \sigma_{\rm s}^\top]$, with $\sigma_{\rm b} \in \R^{n_{\rm b}}$ and $\sigma_{\rm s} \in \R^{n_{\rm s}}$.

In most experiments, the shapes of the electrodes are defined as intersections of the boundary $\partial \Omega$ with a cylinder of radius $R > 0$, oriented such that the normal of $\partial \Omega$ at the center of the electrode lies on the central axis of the cylinder. In the computational model, the boundary of this disk shape is approximated as a polygon on the triangular surface mesh.  This is the electrode model assumed in the rest of this subsection; the required modifications for one of the water tanks that has rectangular electrodes are briefly considered in the next subsection. The contact conductivity $\zeta$ is parametrized similarly to $\sigma$ and $u$, but with the piecewise linear basis functions corresponding to the triangular surface mesh of the tetrahedral FE mesh. 

The contact conductivity function $\zeta$ that we use is a smooth, taking its peak values at the centers of the electrodes and vanishing on the electrode boundaries. We assume $\zeta$ has the same shape $\hat{\zeta}$ on each electrode: 
\begin{equation}
\label{eq:zeta_param}
\zeta|_{E_m}(r_m, \psi_m) =  \zeta_m \hat{\zeta} (r_m, \psi_m) , \qquad r_m \in [0, R), \ \ \psi_m \in [0, 2 \pi),
\end{equation}
where $\zeta_m \in \R_+$ controls the strength of the contact on $E_m$, $m=1, \dots, M$. Here $(r_m, \psi_m)$ characterizes the position on $E_m$ in polar coordinates with respect to the midpoint of $E_m$. The smooth contact shape function is of the form
\begin{equation}
    \label{eq:hatzeta}
    \hat{\zeta} (r, \psi) = \exp\!\left( \tau - \frac{\tau R^p}{R^p - r^p} \right),
\end{equation}
where $\tau, p > 0$ are shape parameters, which take the values $\tau = 6$ and $p = 6$ in our experiments. Like in the case of the internal conductivity, we identify the contact conductivity function with the vector of corresponding peak values $\zeta = [\zeta_1, \dots, \zeta_M]^\top\in \R^M$ in~\eqref{eq:zeta_param}.

These definitions enable discretizing the CEM equations \eqref{eq:cemeqs} and approximately computing related forward solutions by the finite element method (FEM), cf.~\cite{Vauhkonen97}. Additionally, derivatives of measurements with respect to the conductivity parameters~$\sigma_j$ can be computed with the help of \eqref{eq:sderiv} simply by setting $\eta = \varphi_j$. Computing these derivatives for all vertices $j = 1,\dots,n$, we can form the Jacobian matrices required for the projections and reconstructions. When Jacobians with respect to electrode contacts and positions are needed, the derivatives of the measurements with respect to the peak contact values $\zeta_m$, $m=1, \dots, M$, and the spherical coordinates of the electrode centers are computed by resorting to similar techniques detailed in~\cite{Jaaskelainen25}.

\begin{figure}[t]
\center{
  {\includegraphics[width=5.5cm]{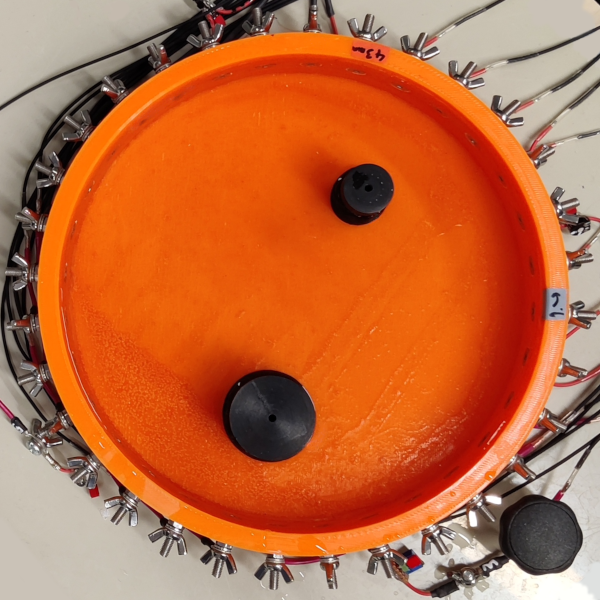}} \qquad
  {\includegraphics[width=5.5cm]{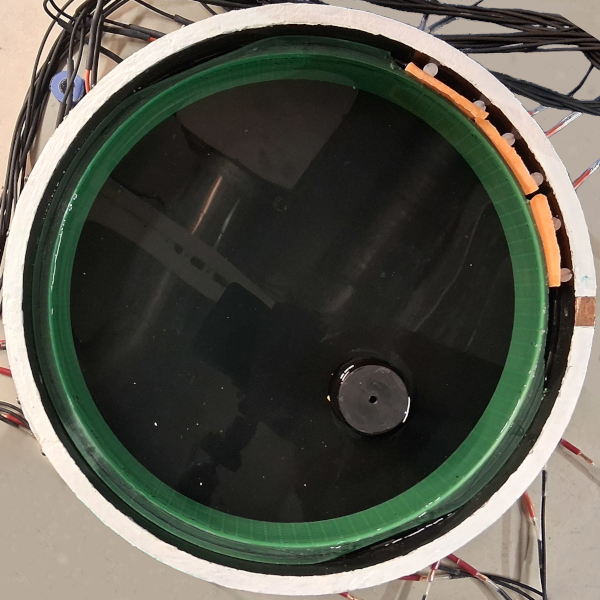}}
  }
\caption{Left: The first water tank with two embedded cylindrical inclusions made of conductive plastic. Right: The second water tank with one embedded cylindrical inclusion made of conductive plastic and sweet potato slices placed between a resistive collar and the boundary of the tank.}
\label{fig:cylindrical_tanks}
\end{figure}

\subsection{Experimental setups}
\label{sec:experimental}
The measurements were performed with the KIT5 stroke EIT device on three different water tanks.  The specifications of the device are given in~\cite{Toivanen21}. In particular, only the electrodes with odd numbers were employed for feeding currents, while all electrodes were used for potential measurements. The  current injections were of 1\,mA at 12\,kHz frequency.

\begin{figure}[t!]
\center{
  {\includegraphics[width=4cm]{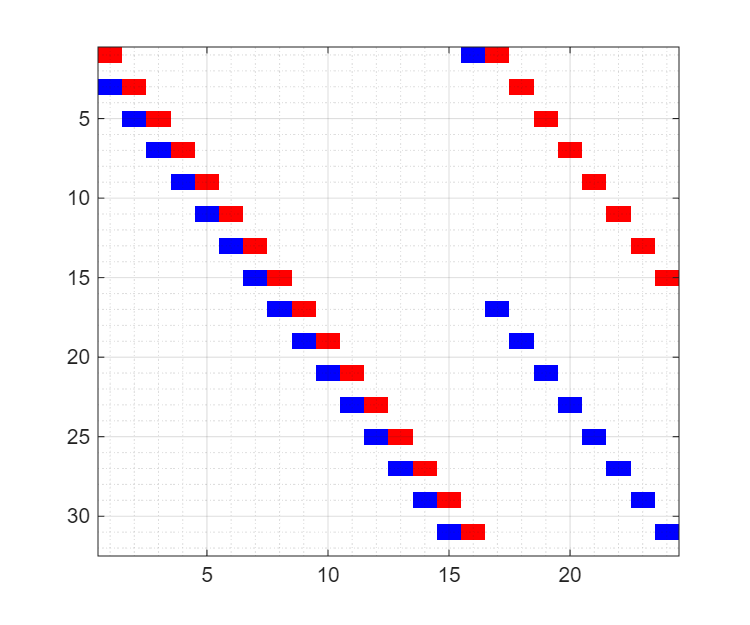}}
  {\includegraphics[width=4cm]{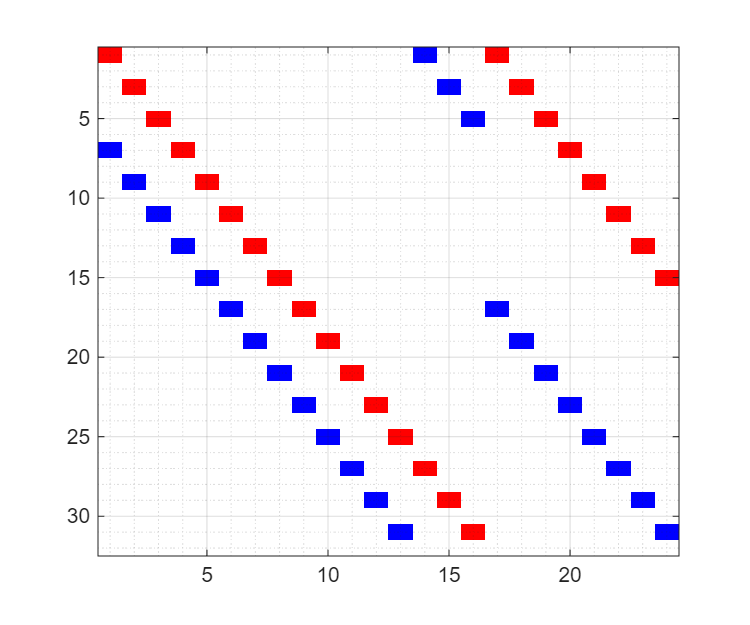}}
  {\includegraphics[width=4cm]{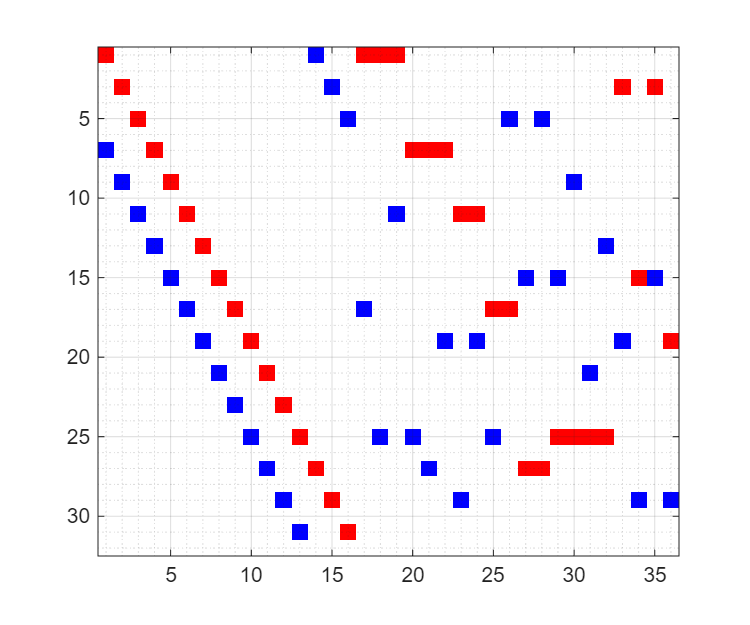}}
  }
\caption{The current patterns used for different water tanks. The vertical axes correspond to electrode numbers and the horizontal axes to different current patterns. Red color indicates the source and blue color the sink of the 1 mA current injection. Left: The first cylindrical tank. Middle: The second cylindrical tank. Right: The head-shaped tank.}
\label{fig:current_patterns}
\end{figure}

The first tank is cylindrical with radius 11.5\,cm, filled with 4.3\,cm of saline solution and with $M=32$ circular electrodes of radius 5\,mm attached equiangularly to its interior surface; see the left image in Figure~\ref{fig:cylindrical_tanks}. The computational domain defined by the water layer is discretized into a FE mesh with roughly $n = 25\,000$ nodes and $120\,000$ tetrahedra, including appropriate refinements of the mesh near the electrode borders to accurately model the contact conductivity function \eqref{eq:hatzeta} and the current flow through the electrodes. The first electrode is the one marked with gray tape in the left image of Figure~\ref{fig:cylindrical_tanks}, and the others are numbered in the counterclockwise direction. The employed current patterns are visualized in the left image of Figure~\ref{fig:current_patterns}.

The second tank is also cylindrical. It also has radius 11.5\,cm, and there are $M=32$ equiangularly positioned rectangular electrodes of width 11\,mm and height 50\,mm on its interior surface; see the right image in Figure~\ref{fig:cylindrical_tanks}. The tank is filled with saline solution up to the top edge of the electrodes. The contact conductivity model \eqref{eq:zeta_param}--\eqref{eq:hatzeta} is modified for the rectangular electrodes as follows: the contact shape $\hat{\zeta}$ is constant in the vertical direction and has the form \eqref{eq:hatzeta} in the horizontal direction around the central line of the considered electrode with $r$ interpreted as the horizontal coordinate and $R$ as the half-width of the electrodes. The computational domain is discretized into a FE mesh with about $n = 11\,000$ nodes and $46\,000$ tetrahedra, with appropriate refinements at the electrodes. The first electrode is the one marked with gray tape in the right image of Figure~\ref{fig:cylindrical_tanks}, and the others are numbered in the counterclockwise direction. The employed current patterns are visualized in the middle image of Figure~\ref{fig:current_patterns}.

\begin{figure}[t]
\center{
{\includegraphics[width=5cm]{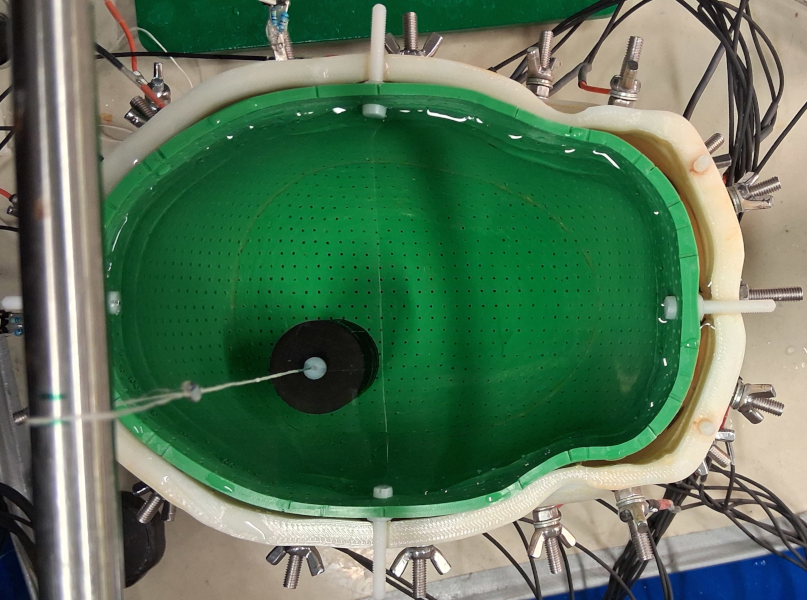}} \qquad
{\includegraphics[width=5cm]{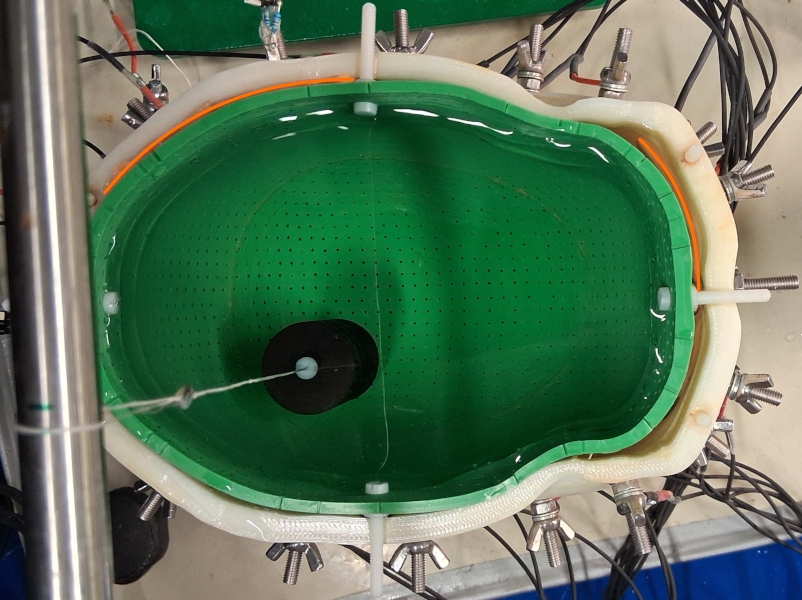}} \linebreak
  {\includegraphics[width=7cm]{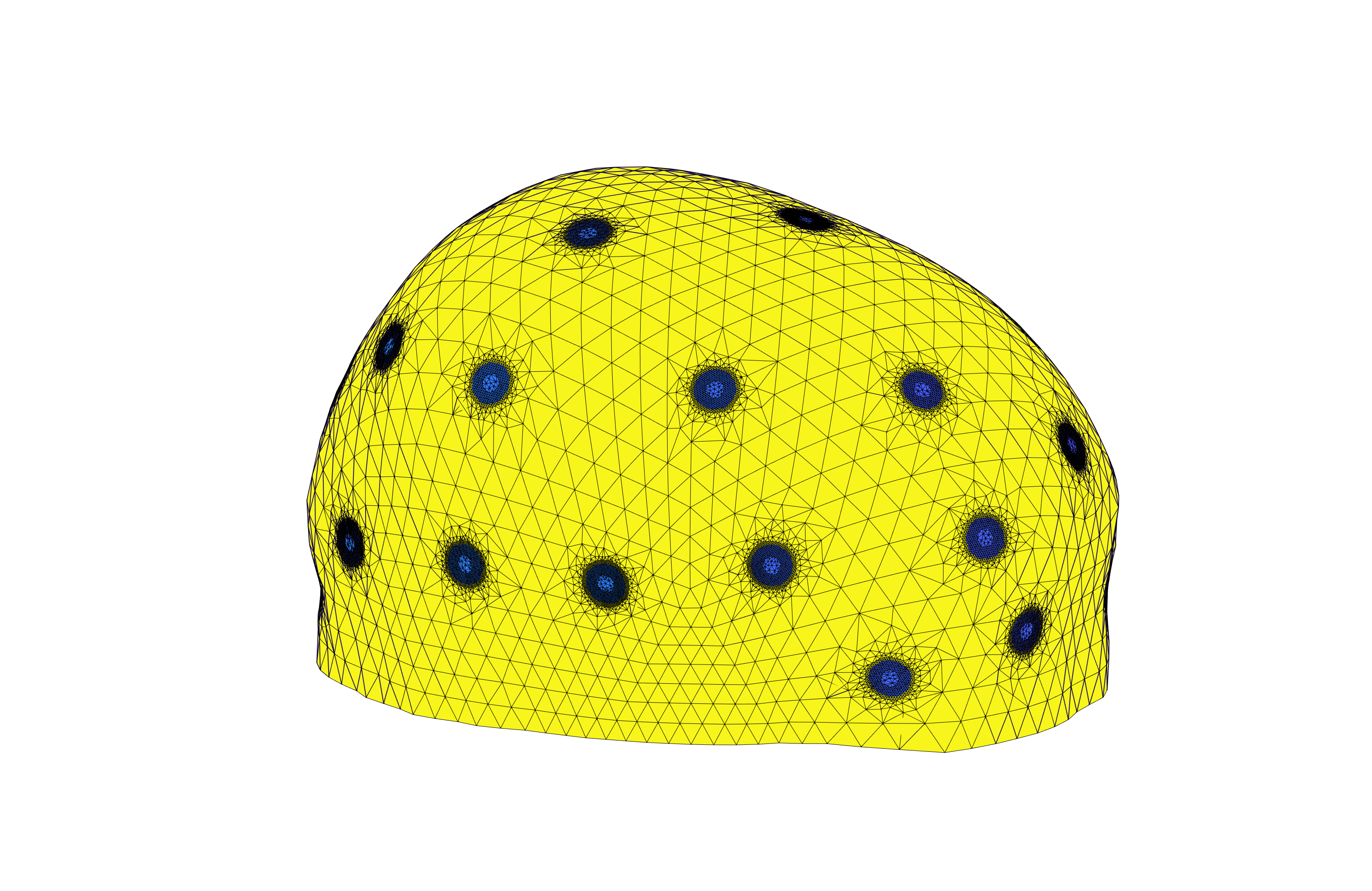}} \hspace{-7mm}
  {\includegraphics[width=5.5cm]{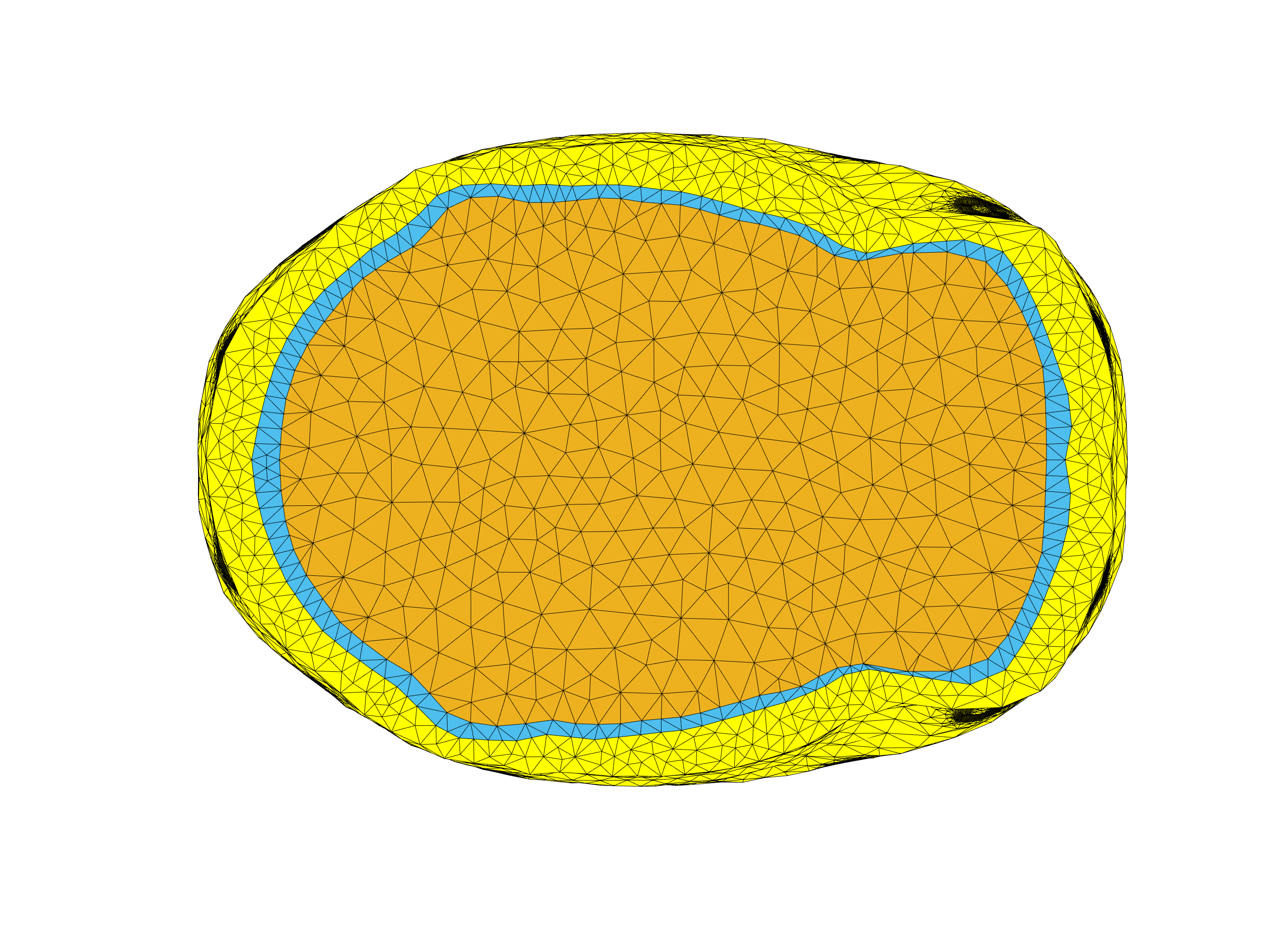}}
  }
\caption{Top left: The head-shaped water tank with a cylindrical embedded inclusion made of conductive plastic inside a 3D-printed skull. Top right: The head-shaped water tank with an embedded inclusion made of conductive plastic and sweet potato slices placed between the skull and the boundary of the tank. Bottom: The FE head model used in the numerical computations.}
\label{fig:headmodel}
\end{figure}

The the third water tank was 3D-printed following the specifications in \cite{Avery17}; see the top row in  Figure~\ref{fig:headmodel}. It mimics the shape of a human head and has the physical dimensions $13.5 \times 19.2 \times 11.9$\,cm. The tank is filled with saline, and a perforated 3D-printed skull is placed inside the tank to imitate the partial shielding effect of the low conductivity human skull. The tank's interior is equipped with $M=32$ electrodes, each of diameter $5$\,mm. A comprehensive description of the tank, including details on the electrode numbering and positioning, is provided in \cite{Avery17}. The FE model for the tank is illustrated on the bottom row of Figure~\ref{fig:headmodel}, and it includes approximately $n=41\,000$ nodes and $180\,000$ tetrahedra. The employed current patterns are visualized in the right image of Figure~\ref{fig:current_patterns}.

\section{Construction of projections}
\label{sec:data}

This section describes how projections can be used to deal with uncertainties arising from unknown nuisance parameters, and how this idea can be utilized to reduce the uncertainty caused by conductivity changes in the RONI, which may,~e.g.,~correspond to the skin layer of a human head. The main principle behind the basic version of our projection method is that the unknown nuisance parameters are assumed to explain every part of the measurement data they are capable of explaining. When estimating the internal conductivity of the imaged object in the ROI, we thus aim to only consider those parts of the available difference data that are definitely caused by deviations of the conductivity in the ROI, and not,~e.g.,~by changes in the contact conductivity or the conductivity of the RONI. 

In a linearized setting, the possible contributions by nuisance parameters are described by the range of the associated Jacobian matrices. In order to eliminate the portion of the data possibly originating from the nuisance parameters, we thus project the data onto the orthogonal complement of the range of the Jacobian --- or alternatively, onto the orthogonal complement of a suitable subspace of the range as explained below. Although the dependence of the measured potentials on the conductivity is nonlinear in EIT, numerical experiments suggest that this approach can regardless perform well even if only rough initial guesses for the unknown parameters are available when computing the needed nuisance Jacobian matrices~\cite{Jaaskelainen25}. 

\subsection{Full projections for low-dimensional parameters}
\label{sec:fullproj}

A projection matrix that aims to eliminate the contribution of a nuisance parameter can be computed as follows. Let $J \in \R^{ML \times K}$ be the Jacobian matrix of the measurement vector $\mathcal{U}$ in \eqref{eq:forward_map} with respect to the considered nuisance parameter evaluated at an initial guess for all unknown parameters. The dimension of this parameter is denoted by $K \in \N$ and assumed to be small,~i.e.,~$K \ll ML$, in this subsection. A projection matrix for projecting onto the orthogonal complement of the range of $J$ can then be formed using the formula
\begin{equation}
\label{eq:proj}
    P = \mathrm{I} - J (J^{\rm \top}\! J)^{-1} J^{\rm \top},
\end{equation}
assuming the nullspace of $J$ is trivial, which is typically the case in EIT if the considered nuisance parameter corresponds to electrode locations or contacts~\cite{Jaaskelainen25}. For low-dimensional nuisance parameters, the evaluation of \eqref{eq:proj} is computationally cheap.

Multiple different parameters can be incorporated into the same projection matrix. With two Jacobian matrices $J_a$ and $J_b$, a matrix $J_{a,b}$ with a range equal to the combined range of the two matrices can be formed by simply concatenating them: 
\begin{equation}
\label{eq:concat}
J_{a,b} = [J_a \ J_b].
\end{equation}
A projection matrix onto the orthogonal complement of the combined range is then formed by setting $J = J_{a,b}$ in \eqref{eq:proj}, assuming that the columns of $J_{a,b}$ are linearly independent. We employed this method in \cite{Jaaskelainen25} to (simultaneously) utilize projections with respect to the contact conductivity parameter $\zeta$ as well as polar $\theta$ and azimuthal $\phi$ angles parameterizing the electrode locations.

\subsection{Prior information based partial projections for the RONI}
\label{sec:skinproj}

In order to employ the projection method to reduce the uncertainty caused by the conductivity of the skin layer or some other RONI, certain modifications are required. Indeed, as any electrode potentials can presumably be (almost) entirely explained by strong enough conductivity changes in a RONI, naively choosing $\sigma_{\rm s}$ as the parameter with respect to which the projection is computed would end up eliminating a majority of the data. To preserve useful information on the conductivity in the~ROI, we are thus forced to perform only a partial projection with respect to the discretized conductivity of the RONI 

A partial orthogonal projection is formed by only considering a subspace of the range of the Jacobian matrix $J_{\sigma_{\rm s}} \in \R^{ML \times n_{\rm s}}$ of the potential measurements with respect to $\sigma_{\rm s} \in \R^{n_{\rm s}}$. The subspace is chosen to include directions that are expected to be affected the most by changes in the conductivity of the RONI. The sensitivity of the measured potentials to a conductivity change in the RONI, say, $\Delta \sigma_{\rm s} \in \R^{n_{\rm s}}$, can be quantified by multiplying it with $J_{\sigma_{\rm s}}$. A simple approach to choosing a basis for the nullspace $\mathcal{N}(P)$ is thus given by selecting orthogonal directions that (sequentially) maximize the norm $\| J_{\sigma_{\rm s} } \Delta \sigma_{\rm s} \|_2$ for $\| \Delta \sigma_{\rm s} \|_2 = 1$. This corresponds to forming the basis for the nullspace out of a chosen number of left singular vectors of $J_{\sigma_{\rm s}}$, or equivalently, eigenvectors of $J_{\sigma_{\rm s}} J_{\sigma_{\rm s}}^{\rm \top}$, corresponding to the largest singular values or eigenvalues, respectively, cf.~\cite{Calvetti25}.

While selecting left singular vectors of $J_{\sigma_{\rm s}}$ as the basis for the nullspace of $P$ provides a method for reducing the dimensionality of the space eliminated by the projection, the restriction $\| \Delta \sigma_{\rm s} \|_2 = 1$ is not particularly meaningful in our setting, where the components of $\Delta \sigma_{\rm s}$ correspond to nodes in the FE mesh, due to the differing sizes of the associated tetrahedra. Furthermore, since we would like to conserve as much of the signal originating from the conductivity changes in the ROI as possible, it is crucial to accurately identify the most important directions to be included in the nullspace of the projection. This can be achieved by incorporating prior information about the conductivity of the RONI into the selection process. A general method to this end is encoding the level of importance, or prior probability, of different conductivity changes $\Delta \sigma_{\rm s}$ in a matrix $B \in \R^{n_{\rm s} \times n_{\rm s}}$. A sensitivity-based weighted projection can then be formed by choosing some number of left singular vectors of $J_{\sigma_{\rm s}} \! B$ corresponding to the largest singular values as the basis of the nullspace $\mathcal{N}(P)$. Denoting $A = B B^{\rm \top}$, the selection of the basis can also be expressed as choosing a symmetric positive definite matrix $A \in \R^{n_{\rm s} \times n_{\rm s}}$ and considering eigenvectors of $J_{\sigma_{\rm s}} \! A J_{\sigma_{\rm s}}^{\rm \top}$. After the orthonormal basis vectors for $\mathcal{N}(P)$, say, $v_1, \dots, v_{K} \in \R^{LM}$, have been selected and stacked as the columns of $V \in \R^{LM \times K}$, the projection itself can be formed via
\[
P = \mathrm{I} - P^\perp, \quad \text{with} \ P^\perp = V V^\top.
\]
Take note that this projection can then be combined with other orthogonal projections, e.g., with respect to the electrode contacts, by employing $P^\perp$ as $J_a$ or $J_b$ in~\eqref{eq:concat}. 

Let us next consider two possible choices for the matrix $A$. For the first case, we assume that any conductivity change of equal $L^2$-norm is equally likely, and thus we wish to eliminate the directions in the data most sensitive to conductivity perturbations that satisfy $\Delta \sigma_{\rm s}^{\rm \top} \! \mathcal{M} \Delta \sigma_{\rm s} = 1$, where $\mathcal{M}$ is the FE mass matrix for the discretized RONI. Since $\mathcal{M}$ is symmetric and positive definite, we may write,~e.g.,~via Cholesky factorization $\mathcal{M} = C^{\rm \top}\! C$, where $C$ is invertible, and $v = C \Delta \sigma_s$ to reduce the problem of finding a basis for $\mathcal{N}(P)$ to the form of deducing orthogonal vectors that sequentially maximize $\| J_{\sigma_{\rm s} } C^{-1} v \|_2$ for $\| v \|_2 = 1$. As previously, the solutions to this problem are given by the first left singular vectors of $J_{\sigma_{\rm s} } C^{-1}$ or, equivalently, the eigenvectors of $J_{\sigma_{\rm s} } \mathcal{M}^{-1} J_{\sigma_{\rm s}}^{\rm \top}$ corresponding to the largest eigenvalues.

Our second approach to using prior information in  constructing the orthogonal projection is assuming a prior probability distribution for $\Delta \sigma_{\rm s}$ with zero mean and a positive definite covariance matrix $\Gamma$. In this case, the linearized contribution of $\Delta \sigma_{\rm s}$ to the measurement data, i.e.~$J_{\sigma_{\rm s}} \Delta \sigma_{\rm s}$, also follows a zero mean distribution with covariance matrix $J_{\sigma_{\rm s} } \Gamma J_{\sigma_{\rm s}}^{\rm \top}$. Finding the orthogonal directions maximizing the expectation of $\| J_{\sigma_{\rm s}} \Delta \sigma_{\rm s} \|_2$ thus corresponds to finding the directions maximizing the variance, which is solved by the eigenvectors of $J_{\sigma_{\rm s} } \Gamma J_{\sigma_{\rm s}}^{\rm \top}$ corresponding to the largest eigenvalues.

\section{Projected reconstruction algorithms}
\label{sec:algorithms}

To judge the effectiveness of projecting away the contribution of conductivity changes in the RONI, we compute reconstructions using data measured from the water tanks described in Section~\ref{sec:experimental}. For this purpose, we use a Bayesian reconstruction algorithm with a smoothened TV prior, which we apply to a linearized version of the SCEM forward problem defined by~\eqref{eq:cemeqs}. This section gives a brief description of the reconstruction algorithm and also explains how orthogonal projections can be incorporated into it.

A Bayesian interpretation of the linearized forward problem leads to
\begin{equation}
\label{eq:lin_mod}
    Y = J_{\sigma_{\rm b}} W + E,
\end{equation}
where the perturbation in the ROI's conductivity $W \in \R^{n_{\rm b}}$, the measured difference potentials~$Y \in \R^{LM}$ and the measurement errors $E\in \R^{LM}$ are random variables. The Jacobian matrix $J_{\sigma_{\rm b}} \in \R^{LM \times n_{\rm b}}$ with respect to $\sigma_{\rm b}$ is evaluated at initial guesses for all unknown parameters. In particular, we compute reconstructions using difference data, that is, $Y$ corresponds to the change in the measured potentials relative to a reference measurement --- depending on the particular experiment, the reference data are either simulated based on initial guesses for the model parameters or measured from the target domain in dynamic difference imaging experiments. We assume the measurement error~$E$ is a zero-mean Gaussian random vector with a positive definite covariance matrix~$\Gamma_E$. The conductivity change $W$ is assumed to follow a smoothened TV prior whose probability density is defined by
\begin{equation}
  \label{eq:TV}
\pi(w) \propto \exp (- \gamma \Psi(w)  ) ,
\end{equation}
where $\gamma > 0$ controls the strength of the prior,
\begin{equation}
\label{eq:TVexp}
\Psi(w) = \int_{\Omega} \upsilon(x) \psi \big(|\nabla w (x)| \big) \, {\rm d} x + \frac{\varepsilon}{2}  |w|^2, \qquad \text{and} \quad \psi(t) = \sqrt{t^2 + T^2} \approx | t |.
\end{equation}
In these equations, $| \, \cdot \, |$ denotes the Euclidean norm/absolute value, $\nabla w$ is to be understood through the representation of $w$ in the FE basis, and $\epsilon, T > 0$ are small parameters such that $\epsilon$ ensures $\pi$ is a proper prior in all considered cases and $T$ guarantees the differentiability of $\psi$. 

The function $\upsilon: \Omega \to \R_+$ in \eqref{eq:TV} is only used in connection to the experiments with the head tank,~i.e.,~we set $\upsilon \equiv 1$ for the circular tanks. For the head tank,
\begin{equation}
    \label{eq:upsilon}
\upsilon(x) = \Big( \frac{1}{2}\big( 1 + \tanh( c_\upsilon  ({\rm dist}(x, \partial \Omega) - b_\upsilon) ) \big) \Big)^{-1}, \qquad x \in \Omega,
\end{equation}
penalizing for conductivity perturbations in the skin and skull layers, which has the effect of ensuring that the algorithm mainly reconstructs conductivity changes in the brain~\cite{Candiani21}. This corresponds to the prior information that the interesting phenomena are expected to occur inside the brain. The parameters $c_\upsilon > 0$ and $b_\upsilon > 0$ are chosen such that $\upsilon$ takes large values in the skin layer and decreases sharply within the skull layer to values close to one; in our experiments with the head tank, the values $c_\upsilon = 300$ and $b_\upsilon = 1$\,cm are used. Take note that the effect of $\upsilon$ on reconstructions is notable only if the skin layer is {\em not} treated as RONI,~i.e.,~when the conductivity change in the skin is also included as a parameter to be reconstructed.

The basic reconstruction algorithm is based on finding the MAP estimate for $W$ by solving for
\begin{equation}
\label{eq:Tikhonovk}
w_{\rm MAP} = {\rm arg} \!\! \! \! \!\! \min_{w \in \R^{n_{\rm b}} \quad} \!\! \! \! \! \frac{1}{2} (y - J_{\sigma_{\rm b}} w )^{\rm \top} \Sigma (y - J_{\sigma_{\rm b}} w ) +  \gamma \, \Psi(w), 
\end{equation}
where $\Sigma = \Gamma_E^{-1}$. When a projection $P$ is utilized in the reconstruction process, we operate with it on~\eqref{eq:lin_mod}, which leads to the projected linearized model
\begin{equation}
\label{eq:proj_lin_mod}
    P Y = P J_{\sigma_{\rm b}} W + P E.
\end{equation}
As the projected measurement noise $P E$ follows a zero-mean Gaussian distribution with covariance matrix $P \Gamma_E P$, the corresponding MAP estimate for $W$ can still be found by solving a minimization problem of the form \eqref{eq:Tikhonovk}, but with the difference that $\Sigma = P \Gamma_E^{-1} P$. This can be seen by computing the pseudoinverse of $P \Gamma_E P$ and applying the properties $P^2 = P = P^{\rm \top}$ of orthogonal projection matrices.

Our reconstructions of the conductivity change in the ROI are computed using the lagged diffusivity method~\cite{Vogel96}, which is an iterative algorithm converging to the MAP estimate of the posterior distribution~\cite{chan1999convergence,dobson1997convergence}. The algorithm can be interpreted to tackle the nonlinearity of the TV prior by successively approximating it with zero-mean Gaussians formed based on the previous iterate~\cite{Helin22}; see also~\cite{Bardsley18,Calvetti08}. The next iterate is then defined as the MAP estimate with the Gaussian approximation as the prior. We present here a brief description of how these iterates are computed, while a more detailed explanation can be found in~\cite{Helin22,Jaaskelainen25}.

The iteration begins from a homogeneous initial guess for the conductivity perturbation $w^{(0)} = 0 \in \R^{n_{\rm b}}$ in the ROI. With the $j$th iterate $w^{(j)}$ in hand, a weight matrix $\Theta(w^{(j)})$ is computed via the formula
\begin{equation}
\label{eq:theta}
    \Theta_{i,j}(w) := \int_{\Omega} \upsilon(x) \frac{\nabla \varphi_i(x) \cdot \nabla \varphi_j(x)}{\sqrt{|\nabla w (x)|^2 + T^2}} \, {\rm d} x + \varepsilon \delta_{i,j} , \qquad i,j=1,\dots ,n_{\rm b}.
\end{equation}
The matrix $\Theta(w^{(j)})^{-1} \in \R^{n_{\rm b} \times n_{\rm b}}$ acts as the covariance matrix for the aforementioned Gaussian approximation for the TV prior at the $(j+1)$th iterate. Successive iterates for $w$ are thus computed using the formula
\begin{equation}
  \label{eq:LD_mean2}
w^{(j+1)} =  \Theta(w^{(j)})^{-1} A^\top \big(\gamma \mathrm{I} + A \Theta(w^{(j)})^{-1} A^\top \big)^{-1} b,
\end{equation}
where the measured difference data $y$ and the Jacobian matrix with respect to the ROI's conductivity define $A$ and $b$ via
\begin{equation}
\label{eq:Amatrix}
A = B J_{\sigma_{\rm b}} \quad \text{and} \quad b = B y.
\end{equation}
Here $B^{\rm \top} \! B = \Sigma$, where, as previously discussed, $\Sigma = \Gamma_E^{-1}$ or $\Sigma = P \Gamma_E^{-1} P$, depending on whether a projection is utilized or not. The iteration \eqref{eq:LD_mean2} is continued by alternating between evaluating equations \eqref{eq:theta} and \eqref{eq:LD_mean2} until satisfactory convergence is observed. In our experiments, very few visual changes between iterates can be seen after ten iterations, which is chosen as the simplistic stopping criterion in our experiments.

\begin{remark}
\label{remark:multilin}
As discussed in~\cite[Remark~5.1]{Jaaskelainen25}, our reconstruction algorithm (with or without projections) can also be applied to the {\em nonlinear} reconstruction problem of EIT via sequential linearizations of the forward model {\em without} recomputing the utilized projection after its initialization. Implementations of such a combination of sequential linearizations and lagged diffusivity iteration can be found in~\cite{Candiani21,Harhanen15}. However, because iteratively accounting for the nonlinearity of the forward model does not lead to significant improvements in the reconstruction quality in our experiments, this paper only considers the simple linearized framework.
\end{remark}

\section{Experiments}
\label{sec:experiments}

In this section, we showcase the results of the experiments we conducted to demonstrate the functionality of our projection method. Three different scenarios are considered, with the experimental setups for each of them described in Section~\ref{sec:experimental}. The first experiment is a basic test, where we project away the effect of the conductivity perturbation in the lower half of the first cylindrical tank shown on the left in Figure~\ref{fig:cylindrical_tanks}. The second experiment studies reducing the contribution of the perturbations near the boundary of the second cylindrical tank on the right in Figure~\ref{fig:cylindrical_tanks}, with the idea of mimicking conductivity changes in the skin layer of a head, but in a far simpler geometry. Finally, the third experiment models skin layer perturbations more accurately in the head-shaped tank portrayed in Figure~\ref{fig:headmodel}.

In all three experiments, we use the algorithm presented in Section \ref{sec:algorithms} to form reconstructions both with and without projections being utilized. In the second and third experiments, where reducing the effect of conductivity perturbations in the skin layer is considered, we use actual difference data between settings with and without embedded inhomogeneities in the tank, while the reconstructions in the first experiment are computed from the difference between actual measurements with inclusions and simulated measurements for an empty tank.
For all three cases, the background conductivity, which is used for linearizing the problem and computing the needed Jacobian matrices, is chosen to correspond to the conductivity of the saline filling the tank. These values were measured to be 0.0215\,S/m, 0.1099\,S/m and 0.0946\,S/m for the first, second and third experiment, respectively. Additionally, in the third experiment we assume a homogeneous background conductivity of 0.03\,S/m in the skull layer, which does not accurately model the conductivity of the perforated 3D-printed skull \cite{Avery17}. The peak values for the contact conductivity functions, needed for computing the Jacobian matrices, are simply assumed to be 500\,${\rm S}/{\rm m}^2$ for all electrodes in all three experiments, which introduces further modeling error to the setting. Take note that the above listed values for the background conductivity and peak contact conductivity values are also used for simulating the reference measurements for the first experiment.

The parameters defining the smoothened TV prior are chosen such that $\gamma = 100$, $T = 10^{-6}$ and $\epsilon$ is the second smallest eigenvalue of $\Theta$, but the precise values for these parameters do not have a significant effect on the results of the experiments. The components of the measurement noise are assumed to be mutually independent with a common standard deviation that is defined to be 0.5\% of the largest difference between two potential values measured from the empty tank for each experiment. 

The partial orthogonal projections with respect to conductivity changes in the RONI are computed according to the second approach described in Section~\ref{sec:skinproj}, incorporating a covariance structure for possible conductivity perturbations. Here we assume the correlation between conductivity values at different points in the RONI depends on the distance between them,
\begin{equation*}
 \Gamma_{i,j} = \varsigma^2 \exp \! \left(-\frac{| x_i - x_j |^2}{2\ell^2} \right), \qquad i,j=n_{\rm b}+1, \dots, n_{\rm b} + n_{\rm s} = n,
\end{equation*}
where $x_i$ denotes a node in the FE mesh. In all experiments, we choose the correlation length as $\ell = 2$\,cm and define the pointwise standard deviation to be $\varsigma = 0.5$\,S/m. The value of $\varsigma$ does not, however, affect the range of the orthogonal projection. The number of orthonormal eigenvectors of $J_{\sigma_{\rm s}} \Gamma J_{\sigma_{\rm s}}^{\top}$ used as the basis for $\mathcal{N}(P)$ is case-dependent and is chosen via trial and error. 

\begin{figure}[ht]
\center{
  {\includegraphics[width=10cm]{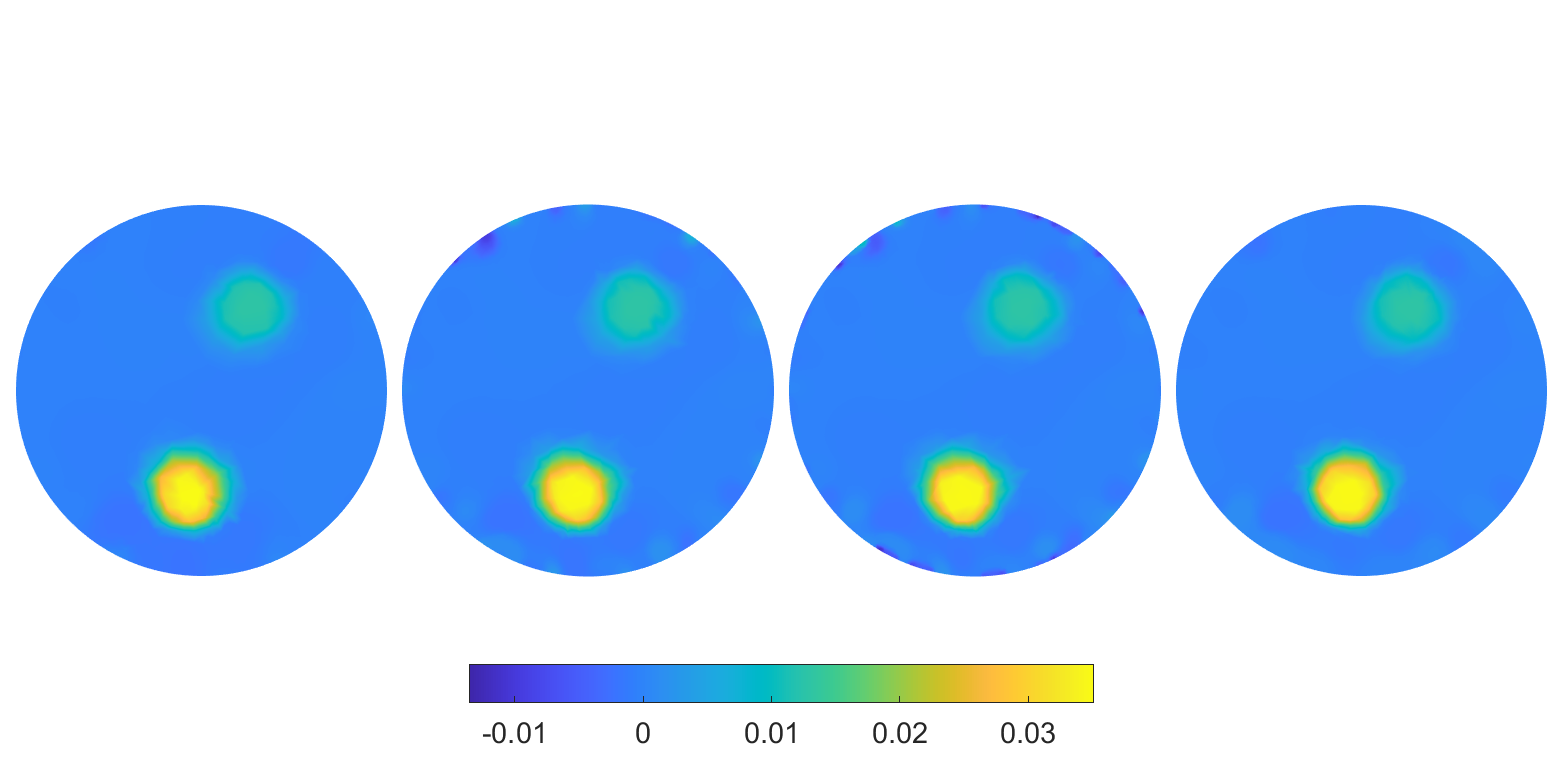}} \\
  \vspace{-2.5mm}
  {\includegraphics[width=10cm]{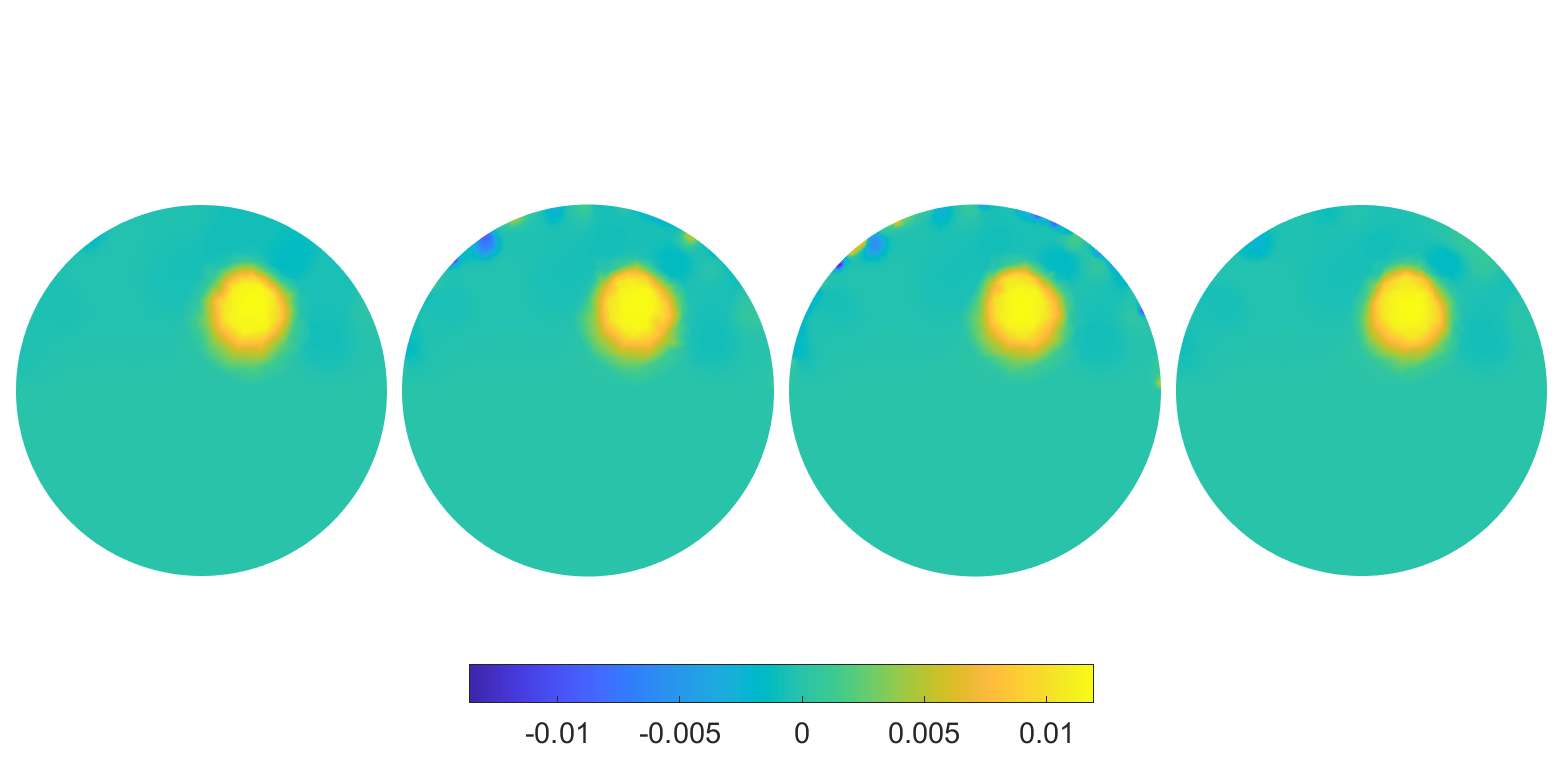}} \\
  \vspace{-2.5mm}
  {\includegraphics[width=10cm]{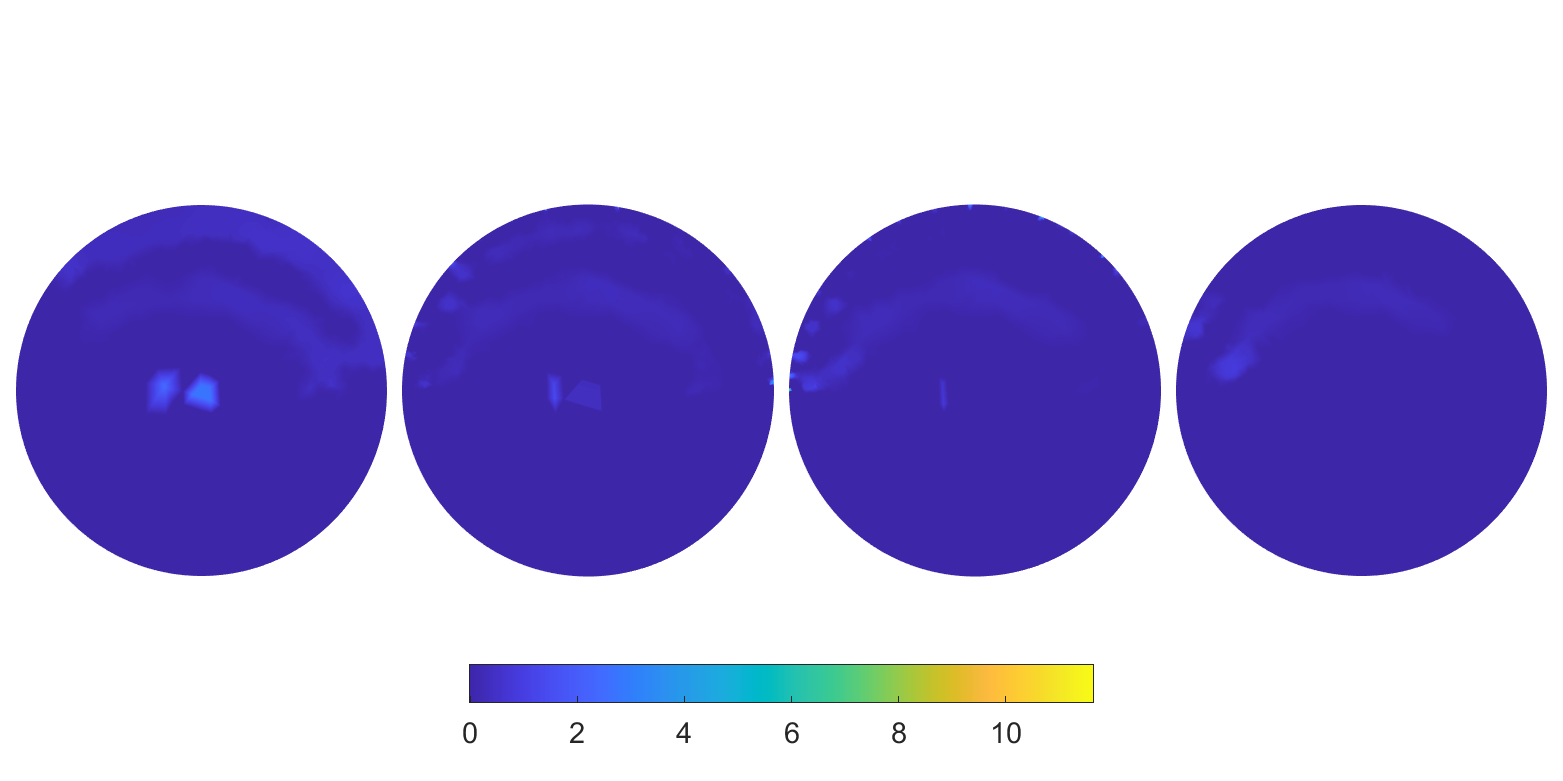}}
  }
\caption{First experiment: Horizontal cross-sections of reconstructions at heights 1\,cm, 2\,cm, 2.5\,cm and 3.5\,cm for the setup in the left image of Figure~\ref{fig:cylindrical_tanks}. The unit of conductivity is S/m. Top: Reconstruction without projections in the whole domain. Middle: Reconstruction when ROI and RONI are, respectively, the top and bottom halves of the tank, and a projection with respect to the RONI is utilized. Bottom: Reconstruction computed only in the top half of the tank without any projections.}
\label{fig:rec_otank}
\end{figure}

\subsection{Experiment~1: First cylindrical water tank}
In the first experiment, we place two distinct cylindrical inclusions of different sizes into a circular water tank, as depicted on the left in Figure~\ref{fig:cylindrical_tanks}. A reconstruction, which is computed assuming the ROI covers the whole of $\Omega$ and no projections are utilized, is visualized on the top row of Figure~\ref{fig:rec_otank}. The reconstruction is of good quality apart from small artifacts close to the domain boundary, caused by the inaccuracy of the employed peak values for the contact conductivity used for simulating the reference measurements. As demonstrated in~\cite{Jaaskelainen25}, employing a projection with respect to the contacts would remove such artifacts, but we choose to not take such an approach in order to better demonstrate the effect of only projecting with respect to the conductivity in the lower half of the tank in what follows.
 
Next, we are interested in reconstructing the top half of the tank only. We achieve this by defining the top half of the tank to be the ROI and the bottom half to be the RONI and employing a projection with respect to the conductivity of the RONI in the algorithm of Section~\ref{sec:algorithms}. When using the 80 most significant eigenvectors of $J_{\sigma_{\rm s}} \Gamma J_{\sigma_{\rm s}}^{\top}$ as the basis for the nullspace of the projection, we obtain the reconstruction shown on the middle row of Figure~\ref{fig:rec_otank}. The larger inclusion is no longer visible in the RONI (as it should not), nor does it cause artifacts in the ROI, while the smaller inclusion appears unaffected, possibly only losing a small amount of its contrast relative to the background. 

As a reference, the bottom row of Figure~\ref{fig:rec_otank} presents a reconstruction that is computed by employing the same division of $\Omega$ into ROI and RONI, but not utilizing a projection to deal with the conductivity perturbation in the RONI. As expected, such a naive approach leads to significant artifacts in the reconstruction of the ROI as the algorithm tries to explain the larger inclusion's contribution to the difference data by conductivity changes in the ROI.

\begin{figure}[t]
\center{
  {\includegraphics[width=12cm]{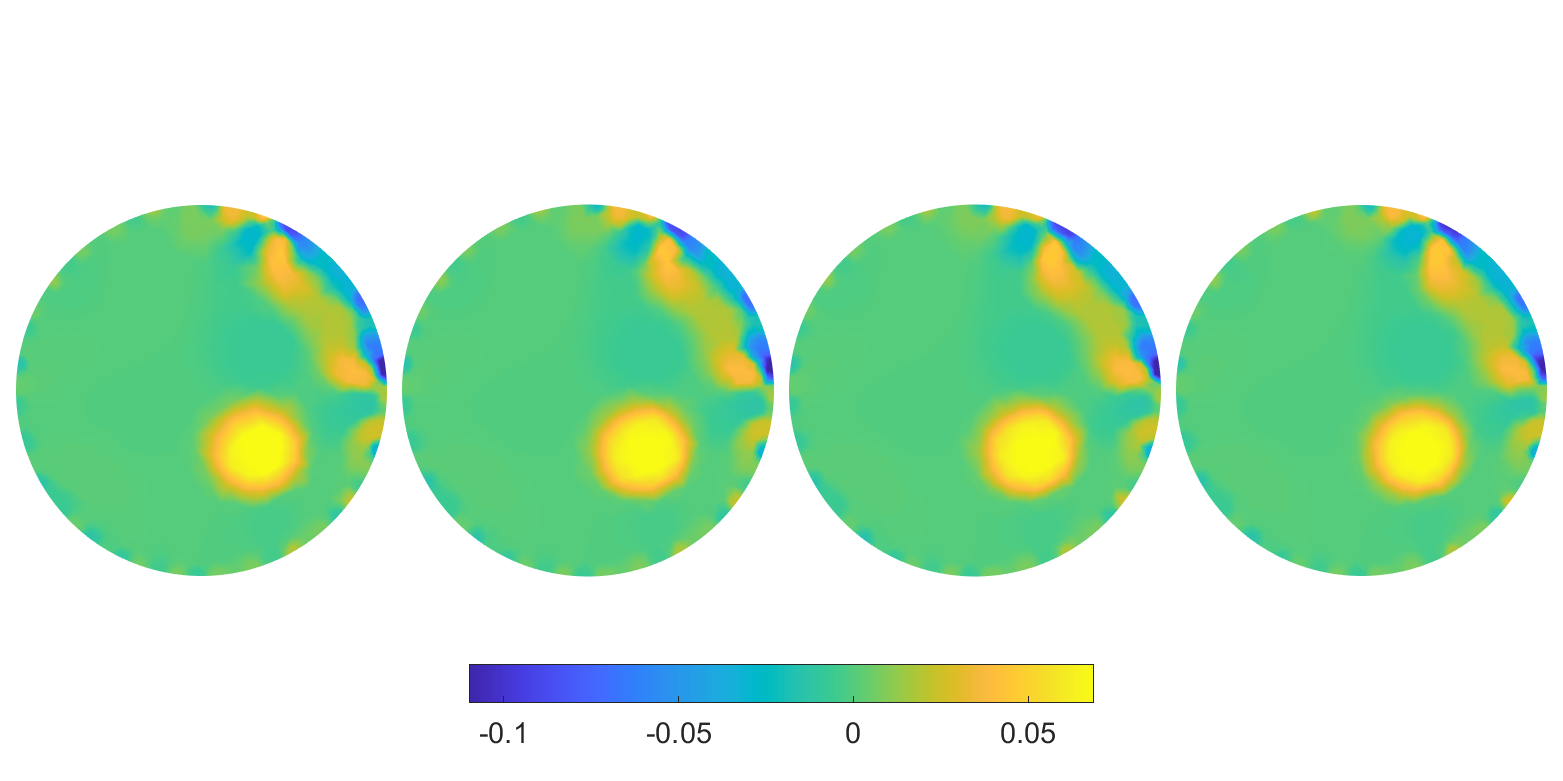}} \\
  \vspace{-3mm}
  {\includegraphics[width=12cm]{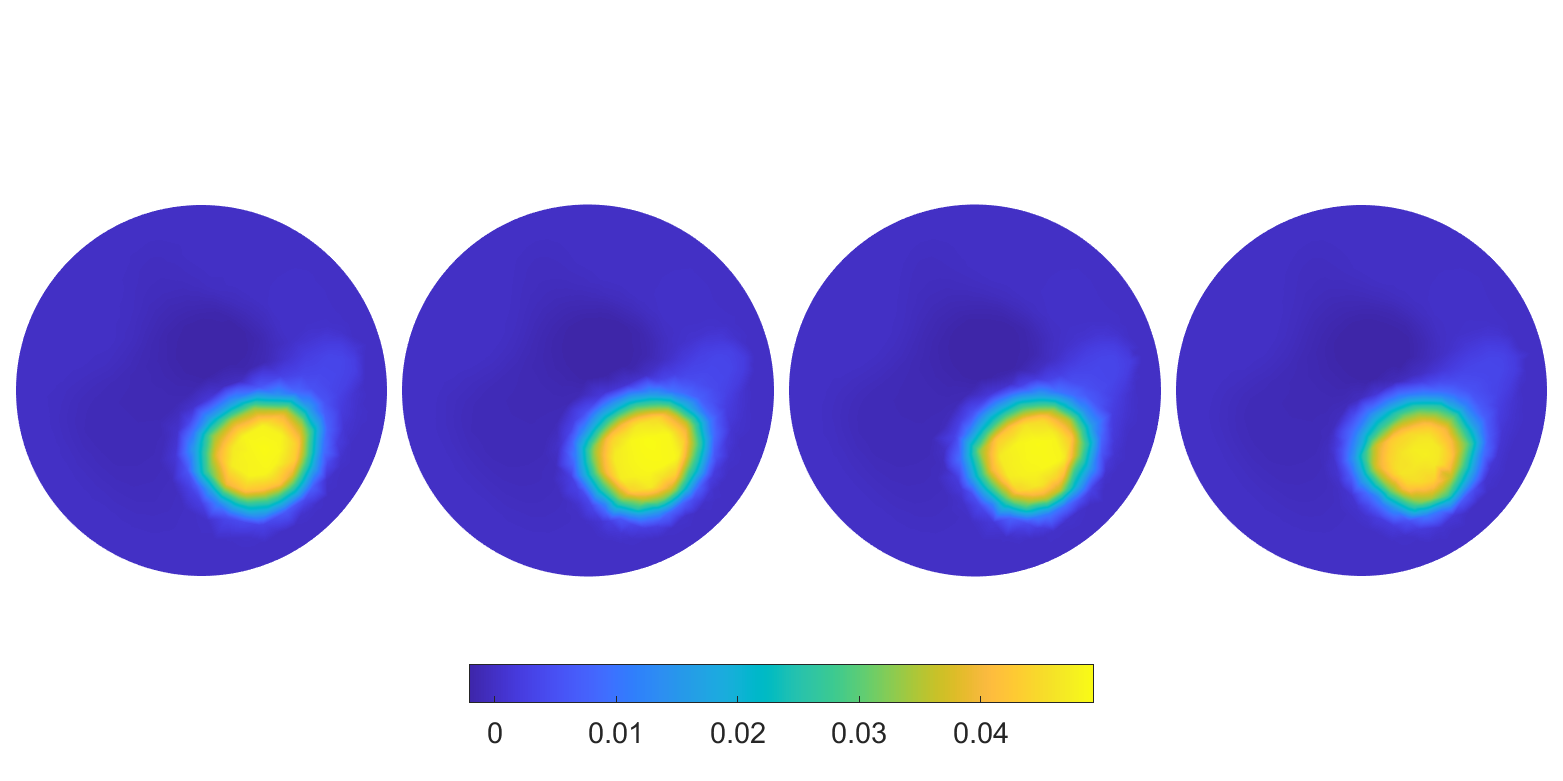}}
  }
\caption{Second experiment: Horizontal cross-sections of reconstructions at heights 1\,cm, 2\,cm, 2.5\,cm and 3.5\,cm for the setup in the right image of Figure~\ref{fig:cylindrical_tanks}. The unit of conductivity is S/m. Top: Reconstruction without projections in the whole domain. Bottom: Reconstruction when ROI is a concentric cylinder of radius 10\,cm, the RONI is its complement, and a projection with respect to the RONI is utilized.}
\label{fig:rec_btank}
\end{figure}

\subsection{Experiment~2: Second cylindrical water tank}
In the second experiment, we take a step towards tackling changes in the conductivity of scalp in stroke monitoring and consider the tank on the right in Figure~\ref{fig:cylindrical_tanks}. There is a single cylindrical conductive inclusion in the tank. In addition, a resistive collar representing the skull is placed in the tank and some sweet potato slices are added onto its exterior, creating conductivity perturbations near the edge of the tank. The reconstructions are computed from difference data between the setting shown in the right image of Figure~\ref{fig:cylindrical_tanks} and a tank that only contains the resistive collar.

Figure \ref{fig:rec_btank} presents the reconstructions from this experiment. The top half of the figure shows the reconstruction when the ROI is $\Omega$ and no projections are used, in which case the cylindrical inclusion is accurately located but the sweet potato slices are also clearly visible close to the tank's boundary. In order to only focus on the conductivity changes inside the collar, we choose the ROI to consist of the points that are at most 10\,cm away from the central vertical axis of the tank and dub the rest of the tank RONI. Utilizing an orthogonal projection $P$ for which the 100 most significant eigenvectors of $J_{\sigma_{\rm s}} \Gamma J_{\sigma_{\rm s}}^{\top}$ form a basis of $\mathcal{N}(P)$, we obtain the reconstruction in the bottom half of Figure \ref{fig:rec_btank}, clearly revealing the location of the cylindrical inclusion without considerable reduction in its contrast. The projection is able to eliminate essentially all uninteresting perturbations caused by the ``skin'', although not as easily as in the first experiment, as demonstrated by the need for a higher dimensional $\mathcal{N}(P)$. The reason for this is presumably the requirement to project away the contribution of a region that borders all the electrodes and thus has a considerable potential impact on the measurements. 

\begin{figure}[t]
\center{
  {\includegraphics[width=12cm]{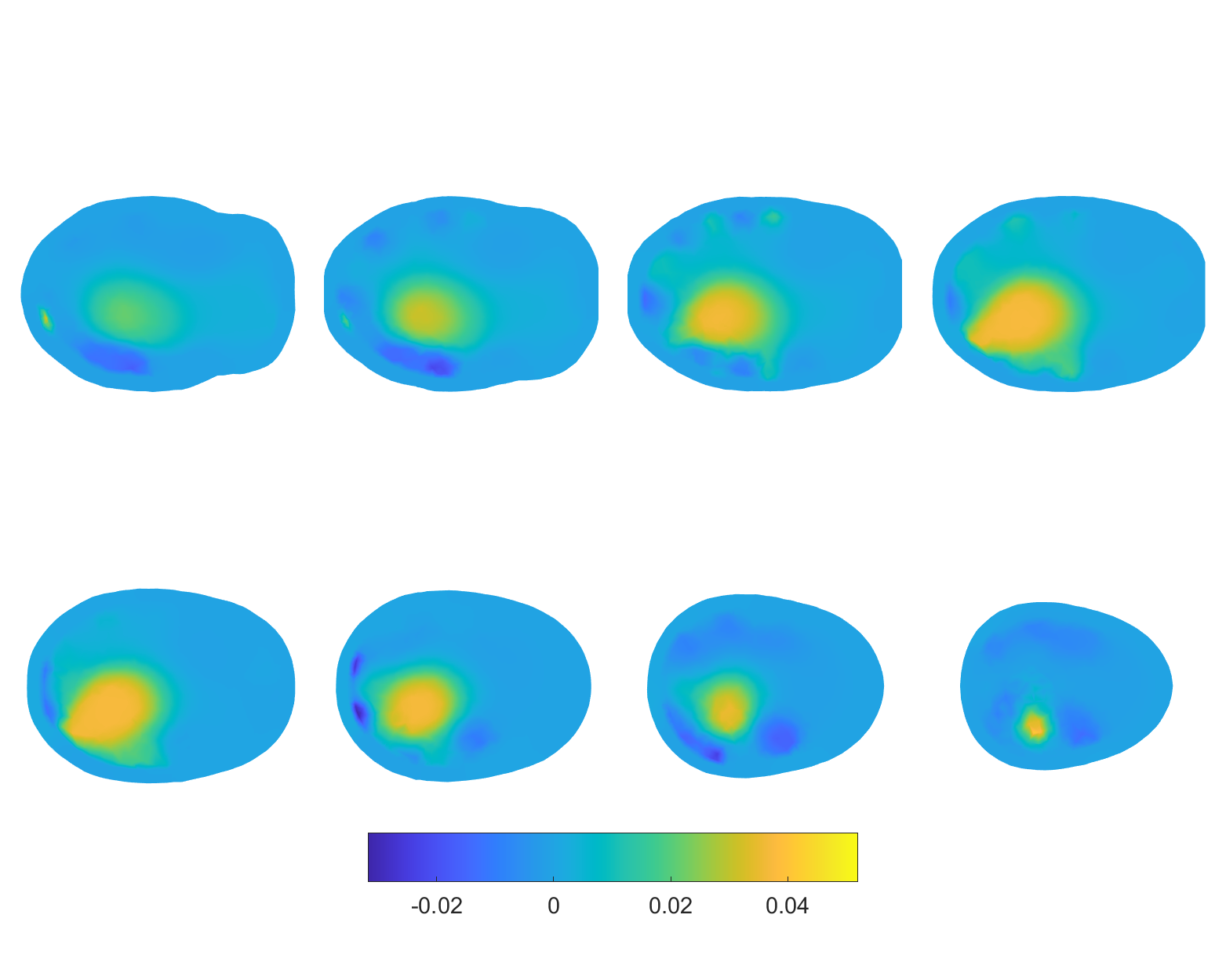}}
  }
\caption{Third experiment without sweet potatoes: Horizontal cross-sections of the reconstruction at heights 2\,cm, 3\,cm, 4\,cm, 5\,cm, 6\,cm, 7\,cm, 8\,cm and 9\,cm for the top left setup in Figure~\ref{fig:headmodel} when no projections are used. The unit of conductivity is~S/m.}
\label{fig:rec_head_nopot_noproj}
\end{figure}

\subsection{Experiment~3: Head-shaped tank}
The third experiment tests projecting away inaccuracies in the measurement geometry and perturbations in the conductivity of the skin layer using measurements from the head-shaped tank illustrated in Figure~\ref{fig:headmodel}. We consider the two setups shown on the top row of Figure~\ref{fig:headmodel}: first, we consider difference data originating from a conductive cylindrical inclusion (left), and subsequently, we tackle a setting where sweet potato slices have also been added between the skull and the exterior boundary of the tank (right). In both cases, the reference measurements correspond to the head tank with the 3D-printed skull in place but no additional inhomogeneities present.

Let us start with the configuration in the top left image of Figure~\ref{fig:headmodel} in case of which we assume that the ROI covers the whole domain and we do not employ any projections with respect to the internal conductivity. Figure \ref{fig:rec_head_nopot_noproj} displays a reconstruction produced by the algorithm of Section~\ref{sec:algorithms} without projections. The highly conductive inclusion is clearly visible, but there are also some additional artifacts that are likely caused by mild errors in the computational model for the head tank and/or the skull placed inside of it. (Recall that the cut-off function included in \eqref{eq:TVexp} forces artifacts closer to the interior of the tank even if they are caused by geometric modeling errors in the exterior boundary shape or electrode positions.)

\begin{figure}[t]
\center{
  {\includegraphics[width=12cm]{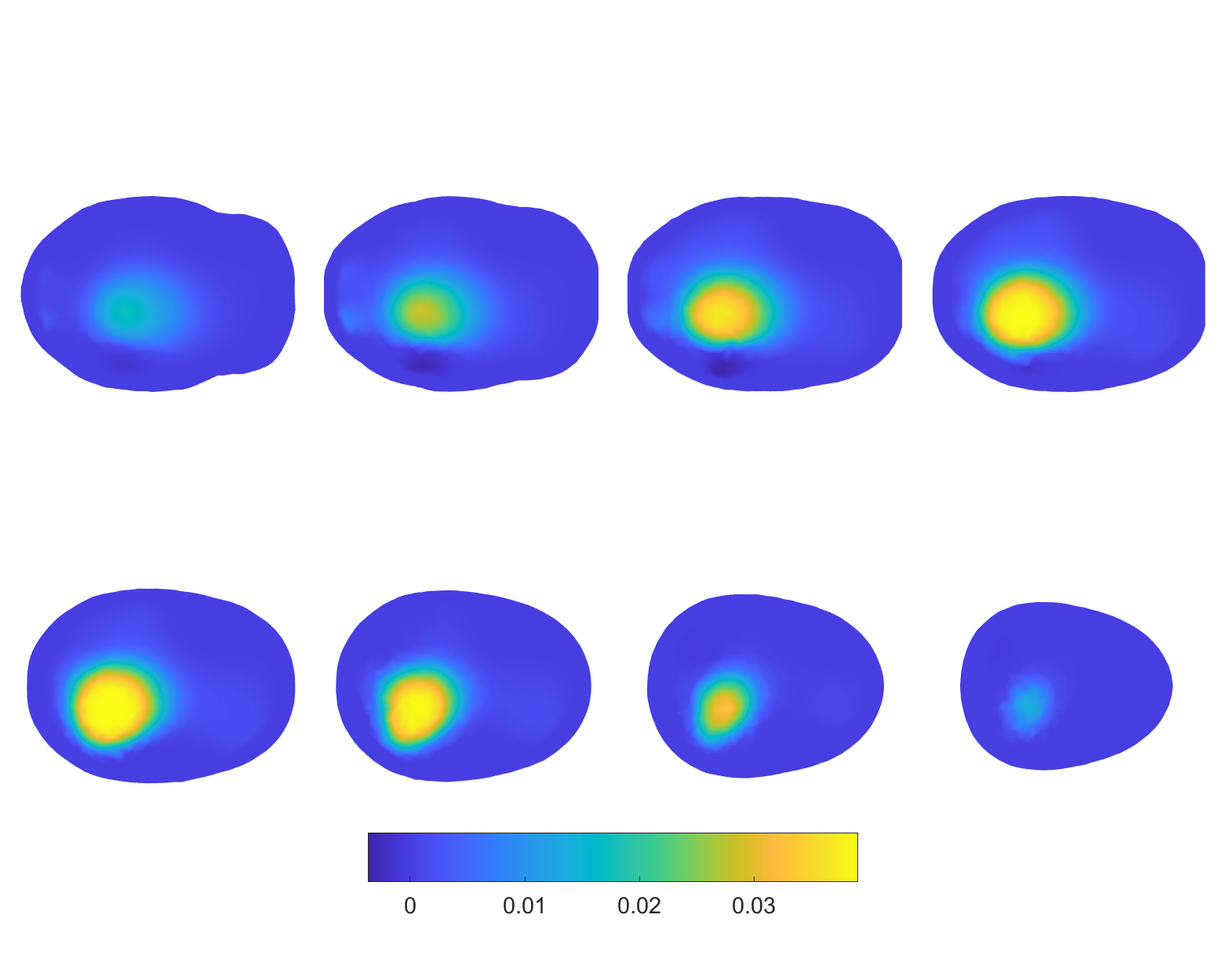}}
  }
\caption{Experiment~3 without sweet potatoes: Horizontal cross-sections of the reconstruction at heights 2\,cm, 3\,cm, 4\,cm, 5\,cm, 6\,cm, 7\,cm, 8\,cm and 9\,cm for the top left setup in Figure~\ref{fig:headmodel} when a projection with respect to the electrode contacts and locations is utilized. The unit of conductivity is S/m.}
\label{fig:rec_head_nopot_zloc}
\end{figure}

Figure \ref{fig:rec_head_nopot_zloc} shows a reconstruction of the same setting when the algorithm exploits projections with respect to the peak values of the contact conductivity and the locations of the electrodes, as introduced in \cite{Jaaskelainen25}. See Section~\ref{sec:fullproj} for the basic ideas of such low-dimensional projections. Although the artifacts in Figure \ref{fig:rec_head_nopot_noproj} were not necessarily caused by modeling errors related to the electrodes, these projections are nevertheless enough to mostly eliminate them and improve the quality of the reconstruction.

\begin{figure}[t]
\center{
  {\includegraphics[width=12cm]{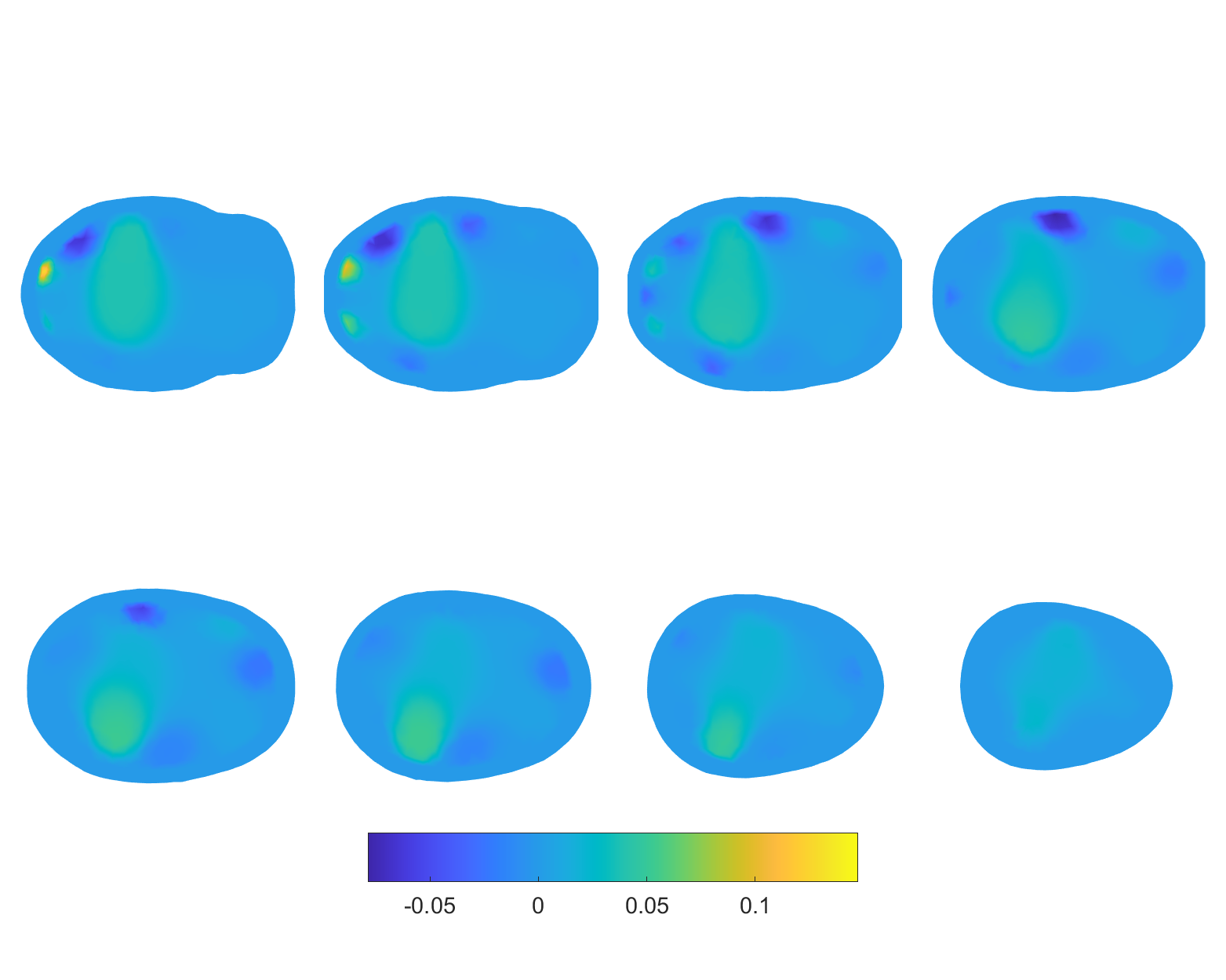}}
  }
\caption{Experiment~3 with sweet potatoes: Horizontal cross-sections of the reconstruction at heights 2\,cm, 3\,cm, 4\,cm, 5\,cm, 6\,cm, 7\,cm, 8\,cm and 9\,cm for the top right setup in Figure~\ref{fig:headmodel} when a projection with respect to the electrode contacts and locations is utilized. The unit of conductivity is S/m.}
\label{fig:rec_head_swtp_zloc}
\end{figure}

Let us then move on to the setup in the top right image of Figure~\ref{fig:headmodel}, with sweet potato slices included between the skull and the boundary of the tank. We start by defining the ROI to be the whole of $\Omega$ and employ orthogonal projections with respect to the electrode contacts and locations in the algorithm. The resulting reconstruction is presented in Figure~\ref{fig:rec_head_swtp_zloc}, indicating that these projections are not able to compensate for all errors caused by the sweet potatoes. Indeed, the reconstruction suffers from severe artifacts along the boundary of the brain and skull layers. Although not documented here, including no projections in the reconstruction algorithm results in even worse reconstruction than the one in Figure~\ref{fig:rec_head_swtp_zloc}.

\begin{figure}[t]
\center{
  {\includegraphics[width=12cm]{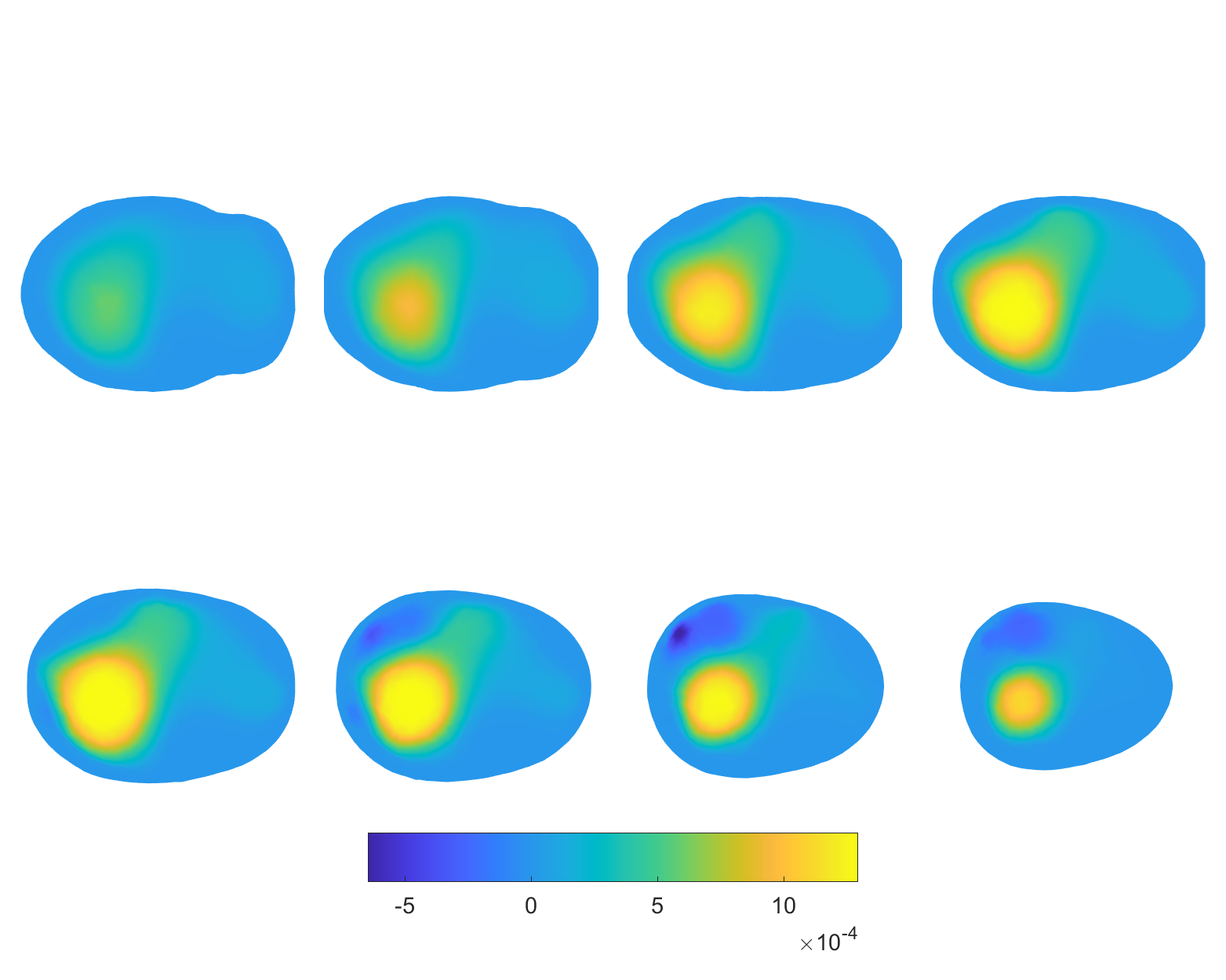}}
  }
\caption{Experiment~3 with sweet potatoes: Horizontal cross-sections of the reconstruction at heights 2\,cm, 3\,cm, 4\,cm, 5\,cm, 6\,cm, 7\,cm, 8\,cm and 9\,cm for the top right setup in Figure~\ref{fig:headmodel} when ROI is the union of the brain and skull, RONI is the skin layer,
and a projection with respect to the RONI is utilized. The unit of conductivity is S/m.}
\label{fig:rec_head_swtp_g180}
\end{figure}

Figure \ref{fig:rec_head_swtp_g180} presents a reconstruction from the same measurements as Figure \ref{fig:rec_head_swtp_zloc}, but instead of electrode parameter projections, we consider the union of the brain and skull layers of the computational model as the ROI and utilize a projection with respect to the conductivity of the RONI,~i.e.,~the skin. Here we use the 180 most significant eigenvectors of $J_{\sigma_{\rm s}} \Gamma J_{\sigma_{\rm s}}^{\top}$ as the basis for $\mathcal{N}(P)$. The resulting reconstruction is considerably better than what was achieved with the electrode parameter projections: the inclusion is now clearly visible and no artifacts can be seen in the background. However, achieving this result required projecting away a much larger portion of the data, which also reduced the available information about the ROI --- note that noiseless measurements with $M=32$ electrodes of which only a half can feed currents correspond to mere 360 degrees of freedom if the self-adjointness of the SCEM forward problem with real conductivities is taken into account. As a consequence, the shape of the reconstructed inclusion has low contrast, and it does not provide a good match with the true inclusion shape.

\section{Concluding remarks}
\label{sec:conclusion}

This work introduced a method for forming localized reconstructions in EIT by extending the ideas introduced in \cite{Jaaskelainen25} for handling unknown contacts at electrodes without estimating their values; see also \cite{Calvetti25}. The method is based on dividing the examined object into a ROI and a RONI and projecting the measured data and the forward map onto the orthogonal complement of a nuisance subspace that is expected to carry most of the information on conductivity changes in the latter. This allows to only compute the reconstruction in the ROI, thus enabling local EIT. We proposed to build the nuisance subspace by combining sensitivity analysis with generic prior information on possible conductivity changes in the RONI via considering left singular vectors of a weighted Jacobian matrix of the measurements with respect to conductivity perturbations in the RONI. 

The functionality of the method was confirmed by numerical experiments based on difference data from three different water tanks, including examples mimicking stroke monitoring by~EIT. In addition to brain imaging, the introduced projection method has potential to serve as an easy-to-implement domain truncation technique in many other applications such as industrial tomography (e.g.,~long pipelines), geological EIT (e.g., half-space domains) and other medical applications (e.g.,~lung imaging).

In our experiments, the dimension of the nuisance subspace was chosen via trial and error, which is a shortcoming that should be tackled in future work; see \cite{Calvetti25} for some ideas in the framework of linear inverse problems. Understanding the interplay of the nuisance subspace with the nonlinearity of EIT and testing the methodology with other imaging modalities also provide interesting avenues for future research.

\bibliographystyle{acm}
\bibliography{skinproj-refs.bib}
\end{document}